\documentclass{amsart}
\usepackage{array}
\usepackage{multirow}
\usepackage{graphicx}
\usepackage{epsfig}
\usepackage{amsfonts,amsmath}
\usepackage{amssymb}
\usepackage{latexsym}
\usepackage{epsf}
\usepackage{graphicx,color,graphics}
\usepackage{amsmath,amssymb}
\usepackage{dsfont}
\usepackage[latin1]{inputenc}  
\newtheorem{theorem}{Theorem}[section]

\newtheorem{lemma}[theorem]{Lemma}

\numberwithin{equation}{section}

\begin{document}

\title[\hfil CAUCHY THERMOELASTIC LAMINATED TIMOSHENKO PROBLEM WITH INTERFACIAL SLIP]{NEW DECAY RATES FOR A CAUCHY THERMOELASTIC LAMINATED TIMOSHENKO PROBLEM WITH INTERFACIAL SLIP UNDER FOURIER OR CATTANEO LAWS}
\author[A. Guesmia]{Aissa Guesmia}
\maketitle

\pagenumbering{arabic}

\begin{center}
Institut Elie Cartan de Lorraine, UMR 7502, Universit\'e de Lorraine\\
3 Rue Augustin Fresnel, BP 45112, 57073 Metz Cedex 03, France
\end{center}
\begin{abstract}
The objective of the present paper is to investigate the decay of solutions for a laminated Timoshenko beam
with interfacial slip in the whole space $\mathbb{R}$ subject to a thermal effect acting only on one component modelled by either Fourier or Cattaneo law. When the thermal effect is acting via the second or third component of the laminated Timoshenko beam (rotation angle displacement or dynamic of the slip), we obtain that both systems, Timoshenko-Fourier and
Timoshenko-Cattaneo systems, satisfy the same polynomial stability estimates in the $L^2$-norm of the solution and its higher order derivatives with respect of the space variable. The decay rate depends on the regularity of the initial data. In addition, the presence and absence of the regularity-loss type property are determined by some relations between the parameters of systems. However, when the thermal effect is acting via the first comoponent of the system (transversal displacement), a new stability condition is introduced for both Timoshenko-Fourier and Timoshenko-Cattaneo systems. This stability condition is in the form of threshold between polynomial stability and convergence to zero. To prove our results, we use the energy method in Fourier space combined with judicious choices of weight functions to build appropriate Lyapunov functionals.
\end{abstract}

\vskip0,1truecm
Energy method, Fourier analysis.
\vskip0,1truecm
{\bf MSC2010:} 34B05, 34D05, 34H05.

\renewcommand{\thefootnote}{}
\footnotetext{E-mail addresse: aissa.guesmia@univ-lorraine.fr. } 

\section{Introduction}

In this paper, we investigate the decay properties of a thermoelastic laminated Timoshenko beam 
with interfacial slip in the whole space $\mathbb{R}$ where the thermal effect is modelled by Fourier law or
Cattaneo law. The first system we consider is the coupling of a laminated Timoshenko system with 
a heat conduction described by Fourier law and given by
\begin{align}\label{s21}
\begin{cases}
\varphi_{tt} -k_1\,(\varphi_{x} +\psi +w)_{x} +\tau_1 \gamma \eta_{x}  = 0, \\
\psi_{tt} - k_2\,\psi_{xx} + k_1 \,(\varphi_{x} +\psi +w) + \tau_2 \gamma \eta_{x} = 0, \\
w_{tt} - k_3\,w_{xx} + k_1 (\varphi_{x} +\psi +w) +\tau_3 \gamma \eta_{x} = 0,\\
\eta_t -k_4 \eta_{xx} +\gamma (\tau_1 \varphi_{xt} +\tau_2 \psi_{xt} +\tau_3 w_{xt})=0,
\end{cases}
\end{align}
and the second system of interest is the coupling between a laminated Timoshenko system with a heat 
conduction described by Cattaneo law and given by
\begin{align}\label{s22}
\begin{cases}
\varphi_{tt} -k_1\,(\varphi_{x} +\psi +w)_{x} +\tau_1 \gamma \eta_{x}  = 0, \\
\psi_{tt} - k_2\,\psi_{xx} + k_1 \,(\varphi_{x} +\psi +w) + \tau_2 \gamma \eta_{x} = 0, \\
w_{tt} - k_3\,w_{xx} + k_1 (\varphi_{x} +\psi +w) +\tau_3 \gamma \eta_{x} = 0,\\
\eta_t +k_4 q_x +\gamma (\tau_1 \varphi_{xt} +\tau_2 \psi_{xt} +\tau_3 w_{xt})=0,\\
q_t +k_5 q +k_4 \eta_x =0,
\end{cases}
\end{align}
where $k_1 ,k_2 ,k_3 ,k_4 ,k_5 , \gamma >0$,
$\varphi =\varphi (x,t)$, $\psi=\psi (x,t)$, $\eta=\eta (x,t)$ and $q=q(x,t)$ denoting the transversal displacement and the rotation angle of the beam, the temperature and the heat flow, respectively, $w=w(x,t)$ is proportional to the amount of slip along the interface, so the third equation in \eqref{s21} and \eqref{s22} describes the dynamics of the slip, $x\in \mathbb{R}$ and $t>0$. 
The termal effect $\gamma\eta_x$ is acting only on one equation of the laminated Timoshenko system, so
\begin{equation}\label{tau123}
(\tau_1 ,\tau_2 ,\tau_3 )\in \{(1,0,0), (0,1,0), (0,0,1)\}.
\end{equation} 
Systems \eqref{s21} and \eqref{s22} are, respectively, subject to the initial conditions
\begin{align}\label{i21}
\begin{cases}
(\varphi ,\psi , w,\eta )(x,0)=(\varphi_0 ,\psi_0 , w_0 ,\eta_0 )(x), \\
(\varphi_t ,\psi_t , w_t )(x,0)=(\varphi_1 ,\psi_1 , w_1 )(x)
\end{cases}
\end{align}
and
\begin{align}\label{i22}
\begin{cases}
(\varphi ,\psi , w,\eta , q)(x,0)=(\varphi_0 ,\psi_0 , w_0 ,\eta_0 ,q_0 )(x), \\
(\varphi_t ,\psi_t , w_t )(x,0)=(\varphi_1 ,\psi_1 , w_1 )(x).
\end{cases}
\end{align} 
The main purpose of this paper is to investigate the capacity of the dissipation, generated by the heat conduction $\gamma\eta_x$ via only one equation of the laminated Timoshenko system, to stabilize \eqref{s21} and \eqref{s22}, and to determine its influence on the decay rate of solutions. We will show that the two cases 
\begin{equation*}
(\tau_1 ,\tau_2 ,\tau_3 )=(1,0,0) \quad\hbox{and}\quad (\tau_1 ,\tau_2 ,\tau_3 )\in \{(0,1,0), (0,0,1)\}
\end{equation*}
are completely different in the following sense: 
\vskip0,1truecm
{\bf Case $(\tau_1 ,\tau_2 ,\tau_3 )=(1,0,0)$}: systems \eqref{s21} and \eqref{s22} are stable if and only if 
\begin{equation}\label{k23}
k_2 \ne k_3,
\end{equation}
and when $k_2 \ne k_3$, the following polynomial stability result holds true for \eqref{s21} and \eqref{s22}: for any $N,\,\ell\in \mathbb{N}^*$ with $\ell\leq N$,
$j\in\{0,\,\ldots,\,N-\ell\}$ and $U_{0}\in H^{N}(\mathbb{R})\cap L^{1}(\mathbb{R})$, there exists $c_0 >0$ such that 
\begin{equation}\label{estim1}
\Vert\partial_{x}^{j}U\Vert_{L^{2}(\mathbb{R})} \leq c_0 \,(1 + t)^{-1/12 - j/6}\,\Vert U_{0} \Vert_{L^{1}(\mathbb{R})} + c_0 \,(1 + t)^{-\ell/2}\,\Vert\partial_{x}^{j+\ell}U_{0} \Vert_{L^{2}(\mathbb{R})},\quad\forall t\in \mathbb{R}_+, 
\end{equation} 
where $\partial_{x}^{j}=\frac{\partial^{j}}{\partial x^j}$, and $U$ and $U_0$ are defined in section 2.
\vskip0,1truecm
{\bf Case $(\tau_1 ,\tau_2 ,\tau_3 )\in \{(0,1,0), (0,0,1)\}$}: when the three speeds of wave propagations of the laminated Timoshenko system are equal; that is
\begin{equation}\label{k123}
k_1 =k_2 = k_3,
\end{equation}
systems \eqref{s21} and \eqref{s22} are stable with the following decay rate: for any $N,\,\ell\in \mathbb{N}^*$ with $\ell\leq N$,
$j\in\{0,\,\ldots,\,N-\ell\}$ and $U_{0}\in H^{N}(\mathbb{R})\cap L^{1}(\mathbb{R})$, there exists $c_0 ,\,{\tilde c}_0 >0$ such that
\begin{equation}\label{estim2}
\Vert\partial_{x}^{j}U\Vert_{L^{2}(\mathbb{R})} \leq c_0 \,(1 + t)^{-1/8 - j/4}\,\Vert U_{0} \Vert_{L^{1}(\mathbb{R})} + c_0 e^{-{\tilde c}_0 t}\,\Vert\partial_{x}^{j+\ell}U_{0} \Vert_{L^{2}(\mathbb{R})},\quad\forall t\in \mathbb{R}_+ . 
\end{equation}
If \eqref{k123} is not satisfied, then the following estimate holds true for \eqref{s21} and \eqref{s22}:
\begin{equation}\label{estim3}
\Vert\partial_{x}^{j}U\Vert_{L^{2}(\mathbb{R})} \leq c_0 \,(1 + t)^{-1/8 - j/4}\,\Vert U_{0} \Vert_{L^{1}(\mathbb{R})} + c_0 \,(1 + t)^{-\ell/4}\,\Vert\partial_{x}^{j+\ell}U_{0} \Vert_{L^{2}(\mathbb{R})},\quad\forall t\in \mathbb{R}_+ . 
\end{equation} 
\vskip0,1truecm
It is well known in the literature that the behavior of the Fourier transform of $U$ in
the low frequency region determines the rate of decay of $U$, while its behavior in the high
frequency rigion imposes a regularity restriction on the initial data known as the regularity-loss property; see 
\cite{4, 8, s4, 6, 10, 7}. It seems that the dissipation generated by the heat conduction is so weak in the high frequency region that it leads to the regularity-loss property in the estimates \eqref{estim1} and \eqref{estim3}. On the other hand, the restriction \eqref{k23} and the fact that the decay rate in \eqref{estim1} is smaller than the one in \eqref{estim2} and \eqref{estim3} indicate that the effect of the heat conduction is better propagated to the whole system from the second or third equation of the laminated Timoshenko system than from the first one.
\vskip0,1truecm
A model describing laminated Timoshenko beams with interfacial slip based on the Timoshenko theory (see, for example, 
\cite{hans, hasp, lota1}) is given by
\begin{align}\label{Ls1}
\begin{cases}
\varphi_{tt} -k_1\,(\varphi_{x} +\psi +w)_{x} = 0, \\
\psi_{tt} - k_2\,\psi_{xx} + k_1 \,(\varphi_{x} +\psi +w) = 0, \\
w_{tt} - k_3\,w_{xx} + k_1 (\varphi_{x} +\psi +w) = 0
\end{cases}
\end{align} 
and can be derived from the following more general model of Bresse-type: 
\begin{align}\label{Bs1}
\begin{cases}
\varphi_{tt} -k_1 \,(\varphi_x +\psi +lw)_{x} -{\tilde l}k_3 (w_x -{\tilde l}\varphi)= 0,\\
\psi_{tt} - k_2\,\psi_{xx} +k_1\,(\varphi_x +\psi +lw) = 0,\\
w_{tt} - k_3\,(w_x -{\tilde l}\varphi)_{x} + l\,k_1\, (\varphi_x +\psi +lw)= 0,
\end{cases}
\end{align}
where $l$ and ${\tilde l}$ are positive constants. The system \eqref{Bs1} coincides with \eqref{Ls1} when $l=1$ and ${\tilde l}=0$. When $w=l={\tilde l}=0$, the system \eqref{Bs1} is reduced to the following Timoshenko-type system: 
\begin{align}\label{Ts1}
\begin{cases}
\varphi_{tt} -k_1 \,(\varphi_x +\psi)_{x} = 0,\\
\psi_{tt} - k_2\,\psi_{xx} +k_1\,(\varphi_x +\psi) = 0.
\end{cases}
\end{align}
\vskip0,1truecm
The well-posedness as well as the stability questions for \eqref{Ls1}, \eqref{Bs1} and \eqref{Ts1} have been the subject of various studies in the literature, where different controls (dampings, memories, heat conduction, ...) and/or boundary conditions (Dirichlet, Neaumann, mixed, ...) have been used to force the solution to converge to zero when time $t$ goes to infinity, and get information on its speed of convergence. 
\vskip0,1truecm
In case of bounded domains, we refer the reader to, for example, \cite{agg, beim, clx, cdflr, gues2, gues4, gues8, gms, lizh, lota1, lota2, lota3, rapo, rvma, wxy} and the refereces therein. 
\vskip0,1truecm
In case of unbounded domains, the stability of \eqref{Bs1} and \eqref{Ts1} has been also treated in the literature for the last few years. In this direction, we mention the papers \cite{s6, s1, s4, s5, s3, s2} (see also the references therein), where some polynomial stability estimates for $L^2$-norm of solutions have been proved using frictional dampings, heat conduction effects or memory controls. In some particular cases, the optimality of the decay rate was also proved. 
\vskip0,1truecm
Our results in the present paper give extensions from the bounded to the unbounded domain case. The proof is based on the energy method combined with the Fourier analysis (by using the transformation in the Fourier space) and well chosen weight functions.
\vskip0,1truecm
The paper is organized as follows. In Section 2, we formulate \eqref{s21} and \eqref{s22} as first order Cauchy systems and give some preliminaries. In Sections 3 and 4, we prove our polynomial stability estimates for \eqref{s21} and \eqref{s22}, respectively. We end our paper by some general comments and other related issues in Section 5.

\section{Formulation of the problems}

To formulate \eqref{s21} and \eqref{s22} in abstract first order systems, we introduce the new variables
\begin{equation}\label{v21}
u=\varphi_t ,\quad y=\psi_t ,\quad \theta =w_t ,\quad v=\varphi_{x} +\psi +w ,\quad z=\psi_x \quad\hbox{and}\quad\phi =w_x .
\end{equation}
Then, the systems \eqref{s21} and \eqref{s22} can be rewritten, respectively, in the forms
\begin{align}\label{e10}
\begin{cases}
v_{t} - u_{x} - y - \theta = 0, \\
u_{t} - k_1 \,v_{x} +\tau_1 \gamma\,\eta_x= 0, \\
z_{t} - y_{x} = 0, \\
y_{t} - k_2\,z_{x} + k_1 v + \tau_2 \gamma\,\eta_x = 0, \\
\phi_{t} - \theta_{x} = 0,  \\
\theta_{t} - k_3\,\phi_{x} + k_1\,v +\tau_3 \gamma \eta_x = 0, \\
\eta_{t} - k_4 \eta_{xx} +\gamma (\tau_1 u_x + \tau_2 y_x +\tau_3 \theta_x)= 0.
\end{cases}
\end{align}
and 
\begin{align}\label{e20}
\begin{cases}
v_{t} - u_{x} - y - \theta = 0, \\
u_{t} - k_1 \,v_{x} +\tau_1 \gamma\,\eta_x= 0, \\
z_{t} - y_{x} = 0, \\
y_{t} - k_2\,z_{x} + k_1 v + \tau_2 \gamma\,\eta_x = 0, \\
\phi_{t} - \theta_{x} = 0,  \\
\theta_{t} - k_3\,\phi_{x} + k_1\,v +\tau_3 \gamma \eta_x = 0, \\
\eta_{t} +k_4 q_x +\gamma (\tau_1 u_x +\tau_2 y_x +\tau_3 \theta_x)= 0,\\
q_t +k_5 q +k_4 \eta_x =0.
\end{cases}
\end{align}
Now, we define the variable $U$ and its initial data $U_0$ by 
\begin{align*}
U = \begin{cases}
(v,\,u,\,z,\,y,\,\phi,\,\theta,\eta)^{T} \quad&\hbox{for \eqref{e10}}, \\
(v,\,u,\,z,\,y,\,\phi,\,\theta ,\eta, q)^{T}  \quad&\hbox{for \eqref{e20}}
\end{cases}\quad\hbox{and}\quad U_0 = \begin{cases}
(v,\,u,\,z,\,y,\,\phi,\,\theta,\eta )^{T} (\cdot ,0)\quad&\hbox{for \eqref{e10}}, \\
(v,\,u,\,z,\,y,\,\phi,\,\theta, \eta,q)^{T} (\cdot,0)\quad&\hbox{for \eqref{e20}} .
\end{cases}
\end{align*}
The systems \eqref{e10} and \eqref{e20} with the initial conditions \eqref{i21} and \eqref{i22}, respectively, are equivalent to 
\begin{align}\label{e11}
\begin{cases}
U_{t} (x,t)+ A_2 U_{xx} (x,t) + A_1 U_x (x,t) + A_0 U (x,t)= 0,\\
U (x, 0)=U_{0} (x),
\end{cases}  
\end{align}
where, for \eqref{e10},  
\begin{equation}\label{e12}
{A}_2 U_{xx} = \left(
\begin{array}{c}
0\\
\\
0 \\
\\
0   \\
\\
0  \\
\\
0\\
\\
0\\
\\
-k_4 \eta_{xx}
\end{array}
\right),\,\,
{A}_1 U_{x} = \left(
\begin{array}{c}
-u_{x} \\
\\
-k_1\,v_{x} +\tau_1 \gamma\eta_x \\
\\
-\,y_{x}   \\
\\
-\,k_2 \,z_{x} +\tau_2 \gamma\eta_x   \\
\\
-\,\theta_{x} \\
\\
-\,k_3\,\phi_{x} +\tau_3 \gamma\eta_x \\
\\
\gamma (\tau_1 u_x +\tau_2 y_x +\tau_3 \theta_x ) 
\end{array}
\right) \,\,\hbox{and}\,\, A_0 U= \left(
\begin{array}{c}
-y - \theta \\
\\
0  \\
\\
0   \\
\\
k_1\,v     \\
\\
0  \\
\\
k_1\,v\\
\\
0 
\end{array}
\right),
\end{equation}
and for \eqref{e20} 
\begin{equation}\label{e13}
{A}_2 = 0,\quad
{A}_1 U_{x} = \left(
\begin{array}{c}
-\,u_{x} \\
\\
-\,k_1 \,v_{x} +\tau_1 \gamma\eta_x\\
\\
-\,y_{x} \\
\\
-\,k_2 \,z_{x} +\tau_2 \gamma\eta_x \\
\\
-\,\theta_{x} \\
\\
-\,k_3 \,\phi_{x} +\tau_3 \gamma\eta_x\\
\\
k_4 q_x + \gamma (\tau_1 u_x +\tau_2 y_x +\tau_3 \theta_x ) \\
\\
k_4 \eta_x 
\end{array}
\right) \quad\hbox{and}\quad A_0 U= \left(
\begin{array}{c}
-\,y - \theta \\
\\
0 \\
\\
0   \\
\\
k_1\,v  \\
\\
0  \\
\\
k_1\,v \\
\\
0\\
\\
k_5 q 
\end{array}
\right).
\end{equation}
For a function $h :\mathbb{R}\to \mathbb{C}$, $Re\, h$, $Im\, h$, $\bar{h}$ and $\widehat{h}$
denote, respectively, the real part of $h$, the imaginary part of $h$, the conjugate of $h$ and the Fourier transformation of $h$. 
Using the Fourier transformation (with respect to the space variable $x$) to transform \eqref{e11} in the Fourier space, we obtain the following first order Cauchy system:
\begin{equation}\label{g6}
\begin{aligned}
\begin{cases}
\widehat{U}_{t} (\xi,\,t)-\xi^2\,{A}_2\widehat{U} (\xi,\,t)+ i\,\xi\,{A}_1 \widehat{U} (\xi,\,t) + A_0 \widehat{U} (\xi,\,t) = 0, 
\quad &\xi\in\mathbb{R},\,\,t>0,\\
\widehat{U}(\xi,\,0) = \widehat{U}_{0}(\xi), \quad &\xi\in\mathbb{R}.
\end{cases}
\end{aligned}
\end{equation}
The solution of \eqref{g6} is given by
\begin{equation}\label{e14}
\widehat{U}(\xi,\,t) = e^{-\,(-\xi^2\,{A}_2 +i\,\xi\,{A}_1 + A_0 )\,t}\ \widehat{U}_{0}(\xi).
\end{equation}
The energy $\widehat{E}$ associated with \eqref{g6} is defined by
\begin{equation}\label{E21}
\widehat{E}(\xi,\,t) = \frac{1}{2}\left[k_1\,\vert\widehat{v}\vert^2 + \vert\widehat{u}\vert^{2} + k_2\,\vert\widehat{z}\vert^{2} + \vert\widehat{y}\vert^{2} + k_3\vert\widehat{\phi}\vert^{2} + \vert\widehat{\theta}\vert^{2} +\vert\widehat{\eta}\vert^{2} \right]
\end{equation}
in case \eqref{e10}, and
\begin{eqnarray}\label{E22}
\widehat{E}(\xi,\,t) =\frac{1}{2}\left[k_1 \,\vert\widehat{v}\vert^2 + \vert\widehat{u}\vert^{2} + k_2 \,\vert\widehat{z}\vert^{2} + \vert\widehat{y}\vert^{2} + k_3\vert\widehat{\phi}\vert^{2} + \vert\widehat{\theta}\vert^{2} +\vert\widehat{\eta}\vert^{2} +\vert\widehat{q}\vert^{2}\right] 
\end{eqnarray}
in case \eqref{e20}. System \eqref{g6} is dissipative, since
\begin{equation}\label{Ep21}
\frac{d}{dt}\widehat{E}(\xi,\,t) = -k_4\xi^2 \vert\widehat{\eta}\vert^{2} 
\end{equation}
in case \eqref{e10}, and
\begin{equation}\label{Ep22}
\frac{d}{dt}\widehat{E}(\xi,\,t) = -k_5 \vert\widehat{q}\vert^{2} 
\end{equation}
in case \eqref{e20}. Indeed, first, we remember the following two trivial identities which will be frequently used in this paper: for any two differentiable functions 
$h ,\,d :\mathbb{R}\to \mathbb{C}$, we have 
\begin{equation}\label{hd1}
\frac{d}{dt} Re\,(h {\bar d}) = Re\,(h_t {\bar d} +d_t {\bar h})
\end{equation}
and 
\begin{equation}\label{hd2}
\frac{d}{dt} Re\,(ih {\bar d}) = Re\,\left[i(h_t {\bar d} -d_t {\bar h})\right].
\end{equation}
In case \eqref{e12}, the first equation in \eqref{g6} is equivalent to 
\begin{align}\label{e101}
\begin{cases}
\widehat{v}_{t} - i\xi\widehat{u} - \widehat{y} - \widehat{\theta} = 0, \\
\widehat{u}_{t} - ik_1 \xi \,\widehat{v} +i\tau_1 \gamma\,\xi\widehat{\eta} = 0, \\
\widehat{z}_{t} - i\xi\widehat{y} = 0, \\
\widehat{y}_{t} - ik_2\xi\,\widehat{z} + k_1 \widehat{v} +i\tau_2 \gamma\,\xi\widehat{\eta} = 0, \\
\widehat{\phi}_{t} - i\xi\widehat{\theta} = 0,  \\
\widehat{\theta}_{t} - ik_3\,\xi\widehat{\phi} + k_1\,\widehat{v} +i\tau_3 \gamma\,\xi\widehat{\eta} = 0,\\
\widehat{\eta}_{t} +k_4 \xi^2 \widehat{\eta} +i\gamma\xi (\tau_1 \widehat{u} +\tau_2 \widehat{y}+\tau_3 \widehat{\theta}).
\end{cases}
\end{align}
Multiplying the equations in \eqref{e101} by $k_1\bar{\widehat{v}}$, $\bar{\widehat{u}}$, $k_2\bar{\widehat{z}}$, $\bar{\widehat{y}}$, $k_3\bar{\widehat{\phi}}$, $\bar{\widehat{\theta}}$ and $\bar{\widehat{\eta}}$, respectively, adding the obtained equations, taking the real part of the resulting expression and using \eqref{hd1}, we obtain \eqref{Ep21}. 
Similarily, in case \eqref{e13}, the first equation of \eqref{g6} is reduced to 
\begin{align}\label{e102}
\begin{cases}
\widehat{v}_{t} - i\xi\widehat{u} - \widehat{y} - \widehat{\theta} = 0, \\
\widehat{u}_{t} - ik_1 \xi \,\widehat{v} +i\tau_1 \gamma\,\xi\widehat{\eta} = 0, \\
\widehat{z}_{t} - i\xi\widehat{y} = 0, \\
\widehat{y}_{t} - i k_2 \xi\,\widehat{z} + k_1 \widehat{v} +i\tau_2 \gamma\,\xi\widehat{\eta} = 0, \\
\widehat{\phi}_{t} - i\xi\widehat{\theta} = 0,  \\
\widehat{\theta}_{t} - ik_3\,\xi\widehat{\phi} + k_1\,\widehat{v} +i\tau_3 \gamma\,\xi\widehat{\eta} = 0,\\
\widehat{\eta}_{t} +ik_4\,\xi\widehat{q} +i\gamma\xi (\tau_1 \widehat{u} +\tau_2 \widehat{y}+\tau_3 \widehat{\theta})=0,\\
\widehat{q}_{t} +k_5 \widehat{q} +ik_4\,\xi\widehat{\eta} =0.
\end{cases}
\end{align} 
Multiplying the equations in \eqref{e102} by $k_1 \bar{\widehat{v}}$, $\bar{\widehat{u}}$, 
$k_2 \bar{\widehat{z}}$, $\bar{\widehat{y}}$, $k_3\bar{\widehat{\phi}}$, $\bar{\widehat{\theta}}$, $\bar{\widehat{\eta}}$ and
$\bar{\widehat{q}}$, respectively, adding the obtained equations, taking the real part of the resulting expression and using \eqref{hd1}, we get \eqref{Ep22}.
\vskip0,1truecm
It is clear that the energy $\widehat{E}$ is equivalent to $\vert\widehat{U}\vert^2$ defined in case \eqref{e101} by
\begin{equation*}
\vert\widehat{U}(\xi,\,t)\vert^2 = \vert\widehat{v}\vert^2 + \vert\widehat{u}\vert^{2} + \vert\widehat{z}\vert^{2} + \vert\widehat{y}\vert^{2} + \vert\widehat{\phi}\vert^{2} + \vert\widehat{\theta}\vert^{2} + \vert\widehat{\eta}\vert^{2},
\end{equation*}
and in case \eqref{e102} by
\begin{equation*}
\vert\widehat{U}(\xi,\,t)\vert^2 = \vert\widehat{v}\vert^2 + \vert\widehat{u}\vert^{2} + \vert\widehat{z}\vert^{2} + \vert\widehat{y}\vert^{2} + \vert\widehat{\phi}\vert^{2} + \vert\widehat{\theta}\vert^{2} + \vert\widehat{\eta}\vert^{2}
+ \vert\widehat{q}\vert^{2} .
\end{equation*}
Since, for $\alpha_1 =\frac{1}{2}\min\{k_1 ,k_2 ,k_3 ,1\}$ and $\alpha_2 =\frac{1}{2}\max\{k_1 ,k_2 ,k_3 ,1\}$, we have 
\begin{equation}\label{equiv}
\alpha_1 \vert\widehat{U}(\xi,\,t)\vert^2 \leq \widehat{E}(\xi,\,t) \leq \alpha_2 \vert\widehat{U}(\xi,\,t)\vert^2 ,\quad \forall \xi \in \mathbb{R}, \,\,\forall t\in \mathbb{R}_+ .
\end{equation} 
\vskip0,1truecm
We finish this section by proving two lemmas, which will be also frequently used in the proof of our stability results.
\vskip0,1truecm
\begin{lemma}
Let $\sigma$, $p$ and $r$ be real numbers such that $\sigma >-1$ and $p,\,r>0$. Then there exists $C_{\sigma,p,r}>0$ such that
\begin{equation}\label{g37+}
\int_{0}^{1} \xi^{\sigma}\,e^{-\,r\,t\,\xi^{p}}\ d\xi \leq C_{\sigma,p,r}\,(1 + t)^{-\,(\sigma + 1)/p} ,\quad\forall \,t\in \mathbb{R}_+ .
\end{equation}
\end{lemma}
\vskip0,1truecm
\begin{proof} For $0\leq t\leq 1$, \eqref{g37+} is evident, for any 
$C_{\sigma,p,r}\geq \frac{2^{(\sigma + 1)/p}}{\sigma +1}$. For $t>1$, we have 
\begin{equation*}
\int_{0}^{1} \xi^{\sigma}\,e^{-\,r\,t\,\xi^{p}}\ d\xi = \int_{0}^{1} \xi^{\sigma +1-p}\,e^{-\,r\,t\,\xi^{p}}\,\xi^{p-1}\, d\xi = \int_{0}^{1} (\xi^{p})^{\left(\sigma + 1-p\right)/p}\,e^{-\,r\,t\,\xi^{p}}\,\xi^{p-1}\ d\xi.
\end{equation*}
Taking $\tau = r\,t\,\xi^{p}$. Then 
\begin{equation*}
\xi^{p} = \frac{\tau}{r\,t} \quad\hbox{and}\quad \xi^{p-1}\ d\xi = \frac{1}{p\,r\,t}\ d\tau.
\end{equation*}
Substituting in the above integral, we find
\begin{equation*}
\int_{0}^{1} (\xi^{p})^{\left(\sigma + 1-p\right)/p}\,e^{-\,r\,t\,\xi^{p}}\,\xi^{p-1}\, d\xi = \int_{0}^{r\,t}\left (\frac{\tau}{r\,t}\right)^{\left(\sigma +1-p\right)/p}\,e^{-\,\tau}\,\frac{1}{p\,r\,t}\, d\tau 
\end{equation*}
\begin{equation*}
\leq \frac{1}{p\,(r\,t)^{\left(\sigma + 1\right)/p}}\int_{0}^{+\infty} \tau^{\left(\sigma + 1-p\right)/p}\,e^{-\,\tau}\ d\tau \leq \frac{2^{(\sigma + 1)/p}}{p\,r^{\left(\sigma + 1\right)/p}}C_{\sigma,p}\,(t + 1)^{-\left(\sigma + 1\right)/p},
\end{equation*}
where 
\begin{equation*}
C_{\sigma,p}=\int_{0}^{+\infty} \tau^{\left(\sigma + 1-p\right)/p}\,e^{-\,\tau}\ d\tau,
\end{equation*}
which is a convergent integral, for any $\sigma > -1$ and $p>0$. This completes the proof of \eqref{g37+} with
\begin{equation*}
C_{\sigma,p,r}=\max\left\{\frac{2^{(\sigma + 1)/p}}{\sigma +1}, \frac{2^{(\sigma + 1)/p}}{p\,r^{\left(\sigma + 1\right)/p}}C_{\sigma,p}\right\} .
\end{equation*} 
\end{proof}
\begin{lemma}
For any positive real numbers $\sigma_1,\, \sigma_{2}$ and $\sigma_{3}$, we have
\begin{equation}\label{ineq}
\sup_{|\xi|\geq 1}\vert\xi\vert^{-\sigma_1}\,e^{-\,\sigma_{2}\,t\,\vert\xi\vert^{-\sigma_3}} \leq \left(1+\sigma_{1}/(\sigma_{2}\sigma_{3})\right)^{\sigma_{1}/\sigma_{3}}
\left(1 + t\right)^{-\,\sigma_{1}/\sigma_3},\quad\forall t\in \mathbb{R}_+ .
\end{equation}
\end{lemma}
\vskip0,1truecm
\begin{proof} 
It is clear that \eqref{ineq} is satisfied for $t=0$. Let $t>0$ and $h(x)=x^{-\sigma_1}\,e^{-\,\sigma_{2}\,t\,x^{-\sigma_3}}$, for $x\geq 1$. Direct and simple computations show that
\begin{equation*}
h^{\prime}(x)=(\sigma_2 \sigma_3 t x^{-\sigma_3} -\sigma_{1})x^{-\sigma_1 -1}\,e^{-\,\sigma_{2}\,t\,x^{-\sigma_3}}.
\end{equation*} 
If $t\geq \sigma_{1}/(\sigma_2 \sigma_3 )$, then 
\begin{eqnarray*}
h(x) &\leq & h(((\sigma_2 \sigma_3 t)/\sigma_1)^{1/\sigma_3 })=((\sigma_2 \sigma_3 )/\sigma_1 )^{-\sigma_1/\sigma_3 } e^{-\sigma_1/\sigma_3} \left(1 + 1/t\right)^{\sigma_{1}/\sigma_3}\left(1 + t\right)^{-\,\sigma_{1}/\sigma_3} \nonumber \\
& \leq & ((\sigma_2 \sigma_3 )/\sigma_1)^{-\sigma_1/\sigma_3 }\left(1 + (\sigma_2 \sigma_3 )/\sigma_1\right)^{\sigma_{1}/\sigma_3} \left(1 + t\right)^{-\,\sigma_{1}/\sigma_3} =\left(1+\sigma_{1}/(\sigma_{2}\sigma_{3})\right)^{\sigma_{1}/\sigma_{3}}
\left(1 + t\right)^{-\,\sigma_{1}/\sigma_3}, 
\end{eqnarray*}
which gives \eqref{ineq} by taking $x=|\xi|$. If $0<t< \sigma_1/(\sigma_2 \sigma_3 )$, then  
\begin{equation*}
h(x)\leq h(1)=e^{-\sigma_2 t} \left(1 + t\right)^{\sigma_{1}/\sigma_3}\left(1 + t\right)^{-\,\sigma_{1}/\sigma_3} \leq
(1+\sigma_1/(\sigma_2 \sigma_3 ))^{\sigma_1/\sigma_3 } \left(1 + t\right)^{-\,\sigma_{1}/\sigma_3}, 
\end{equation*}
which implies \eqref{ineq}, for $x=|\xi|$. 
\end{proof}

\section{Stability: Fourier law \eqref{s21}}

This section is dedicated to the investigation of the asymptotic behavior, when time $t$ goes to infinity, of the solution $U$ of \eqref{e11} in case of Fourier law \eqref{s21}. We will prove \eqref{estim1}, \eqref{estim2} and \eqref{estim3} by showing, first, that $\vert\widehat{U}\vert^{2}$ converges exponentially to zero with respect to time $t$. Then, we prove that the solution \eqref{e14} of \eqref{g6} does not converge to zero when $t$ goes to infinity if 
$(\tau_1, \tau_2 ,\tau_3 )=(1,0,0)$ and $k_2 =k_3$.   
\vskip0,1truecm
In this section and in the next one, $C$ denotes a generic positive constant, and $C_{\varepsilon}$ denotes a generic positive constant depending on some positive constant $\varepsilon$. These generic constants can be different from line to line. Before distinguishing between the three cases \eqref{tau123}, we prove several identities, which will play a crucial role in the proofs.  
\vskip0,1truecm
Multiplying \eqref{e101}$_{4}$ and \eqref{e101}$_{3}$ by $i\,\xi\,\overline{\widehat{z}}$ and $-i\,\xi\,\overline{\widehat{y}}$, 
respectively, adding the resulting equations, taking the real part and using \eqref{hd2}, we obtain
\begin{equation}\label{eq31}
\frac{d}{dt}Re\left(i\,\xi\,\widehat{y}\, \overline{\widehat{z}}\right) = \xi^{2}\left(\vert\widehat{y}\vert^{2} -k_2 \vert\widehat{z}\vert^{2}\right) - k_1\,Re\left(i\,\xi\,\widehat{v}\,\overline{\widehat{z}}\right) +\tau_2\gamma \xi^2\,Re\left(\widehat{\eta}\, \overline{\widehat{z}}\right).
\end{equation}
Multiplying \eqref{e101}$_{2}$ and \eqref{e101}$_{1}$ by $i\,\xi\,\overline{\widehat{v}}$ and $-i\,\xi\,\overline{\widehat{u}}$, 
respectively, adding the resulting equations, taking the real part and using \eqref{hd2}, we find
\begin{equation}\label{eq32}
\frac{d}{dt}Re\left(i\,\xi\,\widehat{u}\, \overline{\widehat{v}}\right) = \xi^{2}\left(\vert\widehat{u}\vert^{2} -k_1 \vert\widehat{v}\vert^{2}\right) - Re\left(i\,\xi\,\widehat{y}\,\overline{\widehat{u}}\right) 
- Re\left(i\,\xi\,\widehat{\theta}\,\overline{\widehat{u}}\right)+\tau_1\gamma\xi^2 \,Re\left(\widehat{\eta}\, \overline{\widehat{v}}\right).
\end{equation}
Multiplying \eqref{e101}$_{6}$ and \eqref{e101}$_{5}$ by $i\,\xi\,\overline{\widehat{\phi}}$ and $-i\,\xi\,\overline{\widehat{\theta}}$, respectively, adding the resulting equations, taking the real part and using \eqref{hd2}, we get
\begin{equation}\label{eq33}
\frac{d}{dt}Re\left(i\,\xi\,\widehat{\theta}\, \overline{\widehat{\phi}}\right) = \xi^{2}\left(\vert\widehat{\theta}\vert^{2} -k_3 \vert\widehat{\phi}\vert^{2}\right) -k_1\, Re\left(i\,\xi\,\widehat{v}\,\overline{\widehat{\phi}}\right) 
+\tau_3\gamma\xi^2 \,Re\left(\widehat{\eta}\, \overline{\widehat{\phi}}\right).
\end{equation}
Multiplying \eqref{e101}$_{6}$ and \eqref{e101}$_{1}$ by $-\xi^2\,\overline{\widehat{v}}$ and $-\xi^2\,\overline{\widehat{\theta}}$, 
respectively, adding the resulting equations, taking the real part and using \eqref{hd1}, we have
\begin{eqnarray}\label{eq34}
\frac{d}{dt}Re\left(-\xi^2\,\widehat{\theta}\, \overline{\widehat{v}}\right) & = & \xi^{2}\left(k_1\,\vert\widehat{v}\vert^{2} - \vert\widehat{\theta}\vert^{2}\right) - \xi^2\,Re\left(i\,\xi\,\widehat{u}\,\overline{\widehat{\theta}}\right) 
- k_3\,\xi^2\,Re\left(i\,\xi\,\widehat{\phi}\,\overline{\widehat{v}}\right)\nonumber \\
& & -\xi^2\,Re\left(\widehat{y}\, \overline{\widehat{\theta}}\right)+\tau_3\gamma\,\xi^2 \,Re\left(i\xi\widehat{\eta}\, \overline{\widehat{v}}\right).
\end{eqnarray}
Multiplying \eqref{e101}$_{4}$ and \eqref{e101}$_{1}$ by $-\xi^2\,\overline{\widehat{v}}$ and $-\xi^2\,\overline{\widehat{y}}$, 
respectively, adding the resulting equations, taking the real part and using \eqref{hd1}, we infer that
\begin{eqnarray}\label{eq35}
\frac{d}{dt}Re\left(-\xi^2\,\widehat{y}\, \overline{\widehat{v}}\right) & = & \xi^{2}\left(k_1\,\vert\widehat{v}\vert^{2} - \vert\widehat{y}\vert^{2}\right) - \xi^2\,Re\left(i\,\xi\,\widehat{u}\,\overline{\widehat{y}}\right) 
- k_2\,\xi^2\,Re\left(i\,\xi\,\widehat{z}\,\overline{\widehat{v}}\right) \nonumber \\
& & -\xi^2\,Re\left(\widehat{\theta}\, \overline{\widehat{y}}\right)+\tau_2\gamma\,\xi^2 \,Re\left(i\xi\widehat{\eta}\, \overline{\widehat{v}}\right).
\end{eqnarray}
Multiplying \eqref{e101}$_{3}$ and \eqref{e101}$_{6}$ by $i\,\xi\,\overline{\widehat{\theta}}$ and $-i\,\xi\,\overline{\widehat{z}}$, respectively, adding the resulting equations, taking the real part and using \eqref{hd2}, we entail
\begin{equation}\label{eq36}
\frac{d}{dt}Re\left(i\,\xi\,\widehat{z}\, \overline{\widehat{\theta}}\right) = -\xi^{2}\, Re\left(\widehat{y}\,\overline{\widehat{\theta}}\right)+k_3\,\xi^{2}\, Re\left(\widehat{\phi}\,\overline{\widehat{z}}\right)
+k_1\,Re\left(i\,\xi\,\widehat{v}\,\overline{\widehat{z}}\right)
-\tau_3\gamma \xi^2\,Re\left(\widehat{\eta}\, \overline{\widehat{z}}\right).
\end{equation}
Mltiplying \eqref{e101}$_{5}$ and \eqref{e101}$_{4}$ by $i\,\xi\,\overline{\widehat{y}}$ and $-i\,\xi\,\overline{\widehat{\phi}}$, respectively, adding the resulting equations, taking the real part and using \eqref{hd2}, we arrive at
\begin{equation}\label{eq37}
\frac{d}{dt}Re\left(i\,\xi\,\widehat{\phi}\, \overline{\widehat{y}}\right) = -\xi^{2}\, Re\left(\widehat{\theta}\,\overline{\widehat{y}}\right)+k_2\,\xi^{2}\, Re\left(\widehat{z}\,\overline{\widehat{\phi}}\right)
+k_1\,Re\left(i\,\xi\,\widehat{v}\,\overline{\widehat{\phi}}\right)
-\tau_2\gamma \xi^2\,Re\left(\widehat{\eta}\, \overline{\widehat{\phi}}\right).
\end{equation}
Multiplying \eqref{e101}$_{2}$ and \eqref{e101}$_{3}$ by $-\,\overline{\widehat{z}}$ and $-\,\overline{\widehat{u}}$, 
respectively, adding the resulting equations, taking the real part and using \eqref{hd1}, it follows that
\begin{equation}\label{eq310}
\frac{d}{dt}Re\left(-\,\widehat{u}\, \overline{\widehat{z}}\right) = -k_1\, Re\left(i\,\xi\,\widehat{v}\,\overline{\widehat{z}}\right) - \,Re\left(i\,\xi\,\widehat{y} \,\overline{\widehat{u}}\right)+\tau_1\gamma \,Re\left(i\xi\widehat{\eta}\, \overline{\widehat{z}}\right).
\end{equation}
Finally, multiplying \eqref{e101}$_{2}$ and \eqref{e101}$_{5}$ by $-\,\overline{\widehat{\phi}}$ and $-\,\overline{\widehat{u}}$, respectively, adding the resulting equations, taking the real part and using \eqref{hd1}, it appears that
\begin{equation}\label{eq312}
\frac{d}{dt}Re\left(-\,\widehat{u}\, \overline{\widehat{\phi}}\right) = -k_1\,Re\left(i\,\xi\,\widehat{v}\,\overline{\widehat{\phi}}\right) - \,Re\left(i\,\xi\,\widehat{\theta} \,\overline{\widehat{u}}\right)+\tau_1\gamma \,Re\left(i\xi\widehat{\eta}\, \overline{\widehat{\phi}}\right).
\end{equation} 

\subsection{Case 1: $(\tau_1 ,\tau_2 ,\tau_3)=(1,0,0)$}

We start by presenting the exponential stability result for \eqref{g6} in the next lemma. 
\vskip0,1truecm
\begin{lemma}\label{lemma32}
Assume that $k_2 \ne k_3 $. Let $\widehat{U}$ be the solution \eqref{e14} of \eqref{g6}. Then there exist $c,\,\widetilde{c}>0$ such that 
\begin{equation}\label{11}
\vert\widehat{U}(\xi,\,t)\vert^{2} \leq \widetilde{c}\,e^{-\,c\,f(\xi)\,t}\,\vert\widehat{U}_0 (\xi)\vert^{2},\quad \forall \,\xi\in \mathbb{R},\quad \forall \,t\in \mathbb{R}_+ ,
\end{equation}
where  
\begin{equation}\label{22}
f(\xi) = \frac{\xi^{6}}{1 + \xi^{2} +\xi^{4} +\xi^{6} +\xi^{8}}.
\end{equation}
\end{lemma}
\vskip0,1truecm
\begin{proof}
Multiplying \eqref{e101}$_{2}$ and \eqref{e101}$_{7}$ by $i\,\xi\,\overline{\widehat{\eta}}$ and $-i\,\xi\,\overline{\widehat{u}}$, respectively, adding the resulting equations, taking the real part and using \eqref{hd2}, we get
\begin{equation}\label{equ1}
\frac{d}{dt}Re\left(i\,\xi\,\widehat{u}\, \overline{\widehat{\eta}}\right) = \gamma\xi^{2}\left(\vert\widehat{\eta}\vert^{2} - \vert\widehat{u}\vert^{2}\right) + k_4\xi^2\,Re\left(i\,\xi\,\widehat{\eta}\,\overline{\widehat{u}}\right) -k_1 \xi^2\,Re\left(\widehat{v}\, \overline{\widehat{\eta}}\right).
\end{equation}
Similarily, multiplying \eqref{e101}$_{6}$ and \eqref{e101}$_{7}$ by $\overline{\widehat{\eta}}$ and $\overline{\widehat{\theta}}$, respectively, adding the resulting equations, taking the real part and using \eqref{hd1}, we find
\begin{equation}\label{equ2}
\frac{d}{dt}Re\left(\widehat{\eta}\, \overline{\widehat{\theta}}\right) = \gamma Re\left(i\,\xi\,\widehat{\theta}\,\overline{\widehat{u}}\right) -k_4 \xi^2\,Re\left(\widehat{\eta}\, \overline{\widehat{\theta}}\right) +k_3\,Re\left(i\xi\widehat{\phi}\, \overline{\widehat{\eta}}\right)-k_1 Re\left(\widehat{v}\, \overline{\widehat{\eta}}\right).
\end{equation}
Also, multiplying \eqref{e101}$_{4}$ and \eqref{e101}$_{7}$ by $\overline{\widehat{\eta}}$ and $\overline{\widehat{y}}$, respectively, adding the resulting equations, taking the real part and using \eqref{hd1}, we obtain
\begin{equation}\label{equ3}
\frac{d}{dt}Re\left(\widehat{\eta}\, \overline{\widehat{y}}\right) = -\gamma Re\left(i\,\xi\,\widehat{u}\,\overline{\widehat{y}}\right) -k_4 \xi^2\,Re\left(\widehat{\eta}\, \overline{\widehat{y}}\right) +k_2\,Re\left(i\xi\widehat{z}\, \overline{\widehat{\eta}}\right)-k_1 Re\left(\widehat{v}\, \overline{\widehat{\eta}} \right).
\end{equation}
We define the functional ${F}_0$ as follows:
\begin{eqnarray}\label{F0}
{F}_0 (\xi,\,t) & = & Re\left[i\,\xi\,\left(\lambda_1\,\widehat{y}\, \overline{\widehat{z}}+ \lambda_3\,\widehat{\theta}\, \overline{\widehat{\phi}} +i\lambda_4 \xi\,\widehat{\theta}\, \overline{\widehat{v}}-\frac{(\lambda_4 +1)k_2}{k_2 -k_3}
\,\widehat{z}\, \overline{\widehat{\theta}} +\frac{(\lambda_4 +1)k_3}{k_2 -k_3}\,\widehat{\phi}\, \overline{\widehat{y}} \right)\right] \nonumber \\ 
& & +Re\left(-\xi^2 \,\widehat{y}\, \overline{\widehat{v}}+i\,\lambda_2\,\xi\widehat{u}\, \overline{\widehat{v}}\right),
\end{eqnarray}
where $\lambda_{1} ,\,\lambda_2 ,\,\lambda_{3}$ and $\lambda_{4}$ are positive constants to be defined later ($F_0$ is well defined since $k_2\ne k_3$). By multiplying \eqref{eq31}-\eqref{eq34}, \eqref{eq36} and \eqref{eq37} by $\lambda_1$, $\lambda_2$, $\lambda_3$, $\lambda_4$, $-\frac{(\lambda_4 +1)k_2}{k_2 -k_3}$ and $\frac{(\lambda_4 +1)k_3}{k_2 -k_3}$, respectively, adding the obtained equations and adding \eqref{eq35}, we deduce that
\begin{eqnarray} \label{g26+}
\frac{d}{dt}{F}_0 (\xi,\,t) &=& - \xi^{2}\left(k_2\lambda_1\,\vert\widehat{z}\vert^{2} 
+k_3\lambda_3 \,\vert\widehat{\phi}\vert^{2} +\left(1-\lambda_1 \right)\,\vert\widehat{y}\vert^{2} +\left(\lambda_4 -\lambda_3 \right)\,\vert\widehat{\theta}\vert^{2} +(k_1\lambda_2 -k_1\lambda_4 -k_1 )\,\vert\widehat{v}\vert^{2}\right)  \nonumber \\ 
&&+I_1 Re (i\xi\widehat{v}\, \overline{\widehat{z}})+I_2 Re (i\xi\widehat{v}\, \overline{\widehat{\phi}})+\xi^2 \left[\lambda_2 \vert\widehat{u}\vert^{2} -Re\left(i\xi\,\left(\lambda_4 \widehat{u}\, \overline{\widehat{\theta}}+ \widehat{u}\, \overline{\widehat{y}}\right)\right)\right] \nonumber \\ 
&& - \lambda_2 Re\left[i\xi\,\left(\widehat{y}\, \overline{\widehat{u}}+ \widehat{\theta}\, \overline{\widehat{u}}
+ i\gamma \xi\widehat{\eta}\, \overline{\widehat{v}}\right)\right],
\end{eqnarray}
where 
\begin{equation*}
I_1 =k_2\xi^2 -k_1\lambda_1 -\frac{(\lambda_4 +1)k_1 k_2}{k_2 -k_3}\quad\hbox{and}\quad I_2 = k_3\lambda_4 \xi^2 -k_1\lambda_3 +\frac{(\lambda_4 +1)k_1 k_3}{k_2 -k_3} .
\end{equation*}
We put
\begin{eqnarray}\label{F}
F_1 (\xi,\,t) =\xi^4 \left[F_0 (\xi,\,t) -\frac{1}{k_1} Re\left(I_1\,\widehat{u}\, \overline{\widehat{z}} +I_2\,\widehat{u}\, \overline{\widehat{\phi}}\right)\right].
\end{eqnarray}
Multiplying \eqref{eq310} and \eqref{eq312} by $\frac{I_1}{k_1}$ and $\frac{I_2}{k_1}$, respectively, adding the obtained equations, adding \eqref{g26+} and multiplying the resulting formula by $\xi^4$, we arrive at
\begin{eqnarray} \label{Fprime1}
\frac{d}{dt}{F}_1 (\xi,\,t) &=& -\xi^{6}\left(k_2\lambda_1\,\vert\widehat{z}\vert^{2} 
+k_3\lambda_3\,\vert\widehat{\phi}\vert^{2} +\left(1-\lambda_1 \right)\,\vert\widehat{y}\vert^{2} +\left(\lambda_4 -\lambda_3 \right)\,\vert\widehat{\theta}\vert^{2} +(k_1\lambda_2 -k_1\lambda_4 -k_1 ) \,\vert\widehat{v}\vert^{2}\right)  \nonumber \\
& & +\lambda_2 \xi^6 \vert\widehat{u}\vert^{2} +\gamma \xi^4 Re\left[i\xi\,\left(\frac{I_1}{k_1}\widehat{\eta}\, \overline{\widehat{z}}+ \frac{I_2}{k_1}\widehat{\eta}\, \overline{\widehat{\phi}}\right)+\lambda_2\xi^2 \widehat{\eta}\, \overline{\widehat{v}}\right]+\xi^{4} Re\left(iI_3 \xi\widehat{u}\, \overline{\widehat{\theta}} +iI_4\xi\widehat{u}\, \overline{\widehat{y}}\right),
\end{eqnarray}
where 
\begin{equation*}
I_3 =\lambda_2 +\frac{1}{k_1}I_2 -\lambda_4\xi^2\quad\hbox{and}\quad I_4 = \lambda_2 +\frac{1}{k_1}I_1 -\xi^2.
\end{equation*}
Let $\lambda_5 >0$ and
\begin{eqnarray}\label{FF1}
F (\xi,\,t) =F_1 (\xi,\,t)+\lambda_5 \xi^4 Re\left(i\xi\widehat{u}\, \overline{\widehat{\eta}}\right)+\frac{1}{\gamma}I_3 \xi^{4}  Re\left(\widehat{\eta}\, \overline{\widehat{\theta}}\right)
+\frac{1}{\gamma}I_4 \xi^{4} Re\left(\widehat{\eta}\, \overline{\widehat{y}}\right).
\end{eqnarray}
Multiplying \eqref{equ1}, \eqref{equ2} and \eqref{equ3} by $\lambda_5\xi^4$, $\frac{1}{\gamma}I_3\xi^{4} $ and $\frac{1}{\gamma}I_4 \xi^{4}$, respectively, adding the obtained equations and adding \eqref{Fprime1}, we see that
\begin{eqnarray} \label{FFprime1}
\frac{d}{dt}{F} (\xi,\,t) &=& -\xi^{6}\left(k_2\lambda_1\,\vert\widehat{z}\vert^{2} 
+k_3\lambda_3\,\vert\widehat{\phi}\vert^{2} +\left(1-\lambda_1 \right)\,\vert\widehat{y}\vert^{2} +\left(\lambda_4 -\lambda_3 \right)\,\vert\widehat{\theta}\vert^{2} +(k_1\lambda_2 -k_1\lambda_4 -k_1 ) \,\vert\widehat{v}\vert^{2}\right)  \nonumber \\
& & -(\gamma \lambda_5 -\lambda_2 )\xi^6 \vert\widehat{u}\vert^{2} +\gamma \lambda_5 \xi^6 \vert\widehat{\eta}\vert^{2}
+\lambda_5 \xi^6 Re\left(ik_4 \xi\widehat{\eta}\, \overline{\widehat{u}}-k_1\widehat{v}\, \overline{\widehat{\eta}}\right)\nonumber \\
& & +\frac{1}{\gamma}I_3 \xi^{4} Re\left(ik_3\xi\,\widehat{\phi}\, \overline{\widehat{\eta}}-k_4\xi^2\widehat{\eta}\, \overline{\widehat{\theta}} -k_1 \widehat{v}\, \overline{\widehat{\eta}}\right)
+\frac{1}{\gamma}I_4 \xi^{4} Re\left(ik_2\xi\,\widehat{z}\, \overline{\widehat{\eta}}-k_4\xi^2\widehat{\eta}\, \overline{\widehat{y}} -k_1 \widehat{v}\, \overline{\widehat{\eta}}\right) \nonumber \\
& & +\gamma \xi^4 Re\left[i\xi\,\left(\frac{I_1}{k_1}\widehat{\eta}\, \overline{\widehat{z}}+ \frac{I_2}{k_1}\widehat{\eta}\, \overline{\widehat{\phi}}\right)+\lambda_2\xi^2 \widehat{\eta}\, \overline{\widehat{v}}\right].
\end{eqnarray}
Applying Young's inequality for the terms depending on $\widehat{\eta}$ in \eqref{FFprime1}, it follows that, for any $\varepsilon_0 >0$,
\begin{eqnarray} \label{Fprime2}
\frac{d}{dt}{F} (\xi,\,t) &\leq & -(k_2\lambda_1 -\varepsilon_0 ) \xi^{6}\,\vert\widehat{z}\vert^{2} 
-(k_3\lambda_3 -\varepsilon_0 )\xi^{6}\,\vert\widehat{\phi}\vert^{2} -\left(1-\lambda_1 -\varepsilon_0 \right)\xi^{6}\,\vert\widehat{y}\vert^{2} \nonumber \\ 
& & -\left(\lambda_4 -\lambda_3 -\varepsilon_0 \right)\xi^{6}\,\vert\widehat{\theta}\vert^{2} -(k_1\lambda_2 -k_1\lambda_4 -k_1  -\varepsilon_0 )\xi^{6} \,\vert\widehat{v}\vert^{2} -\left(\gamma \lambda_5 -\lambda_2 -\varepsilon_0 \right)\xi^{6}\,\vert\widehat{u}\vert^{2} \nonumber \\
& & +C_{\varepsilon_0 ,\lambda_1 ,\cdots,\lambda_4 }(1+\xi^2 +\xi^4 +\xi^6 +\xi^8 ) \xi^2\vert\widehat{\eta}\vert^{2} .
\end{eqnarray}
We choose $0<\lambda_{1}<1$, $\lambda_{2} >1$,  $\lambda_{5} >\frac{1}{\gamma}\lambda_{2}$, $0<\lambda_{3}<\lambda_{4}<\lambda_{2} -1$ and 
\begin{equation*}
0<\varepsilon_{0} <\min\left\{k_2\lambda_1 ,k_3\lambda_3 ,1- \lambda_1 ,\lambda_4 - \lambda_3 , k_1\lambda_2 -k_1\lambda_4 -k_1 , \gamma \lambda_5 -\lambda_2\right\}. 
\end{equation*} 
Hence, using the definition \eqref{E21} of $\widehat{E}$, \eqref{Fprime2} leads to, for some positive constant $c_1$,
\begin{eqnarray}\label{g27+} 
\frac{d}{dt}{F}(\xi,\,t) \leq -c_1 \xi^{6}\widehat{E} (\xi,\,t)+ C \,\left(1 +\xi^2 +\xi^{4} + \xi^{6} +\xi^{8} \right)\xi^2\vert\widehat{\eta}\vert^{2} .
\end{eqnarray}
Now, we introduce the {\it Perturbed Energy} ${L}$ as follows:
\begin{equation}\label{g27}
{L}(\xi,\,t) = \lambda\,\widehat{E}(\xi,\,t) + \frac{1}{1 + \xi^{2} + \xi^{4} + \xi^{6} + \xi^{8}}\,{F}(\xi,\,t) ,
\end{equation}
where $\lambda$ is a positive constant to be fixed later. Then from \eqref{Ep21}, \eqref{g27+} and \eqref{g27} we have 
\begin{eqnarray}\label{g30}
\frac{d}{dt}{L}(\xi,\,t) \leq -c_1 f(\xi)\widehat{E}(\xi,\,t)-\left(k_4\,\lambda - C\right)\xi^2\vert\widehat{\eta}\vert^{2} , 
\end{eqnarray}
where $f$ is defined in \eqref{22}. Moreover, using the definitions \eqref{E21}, \eqref{FF1} and \eqref{g27} of $\widehat{E} ,\,{F}$ and ${L}$, respectively, we get, for some $c_{2} >0$ (not depending on $\lambda$),
\begin{equation}\label{l1e} 
\vert{L}(\xi,\,t) - \lambda\,\widehat{E}(\xi,\,t)\vert\leq \frac{c_{2}\,(\xi^4 +\vert\xi\vert^5 +\xi^6 )}{1 + \xi^{2} +\xi^{4} +\xi^{6}+\xi^{8}}\,\widehat{E}(\xi,\,t) \leq 3c_{2}\,\widehat{E}(\xi,\,t).
\end{equation} 
Therefore, for $\lambda$ large enough so that $\lambda >\max\left\{\frac{C}{k_4},\,3c_{2}\right\}$, we deduce from \eqref{g30} and \eqref{l1e} that 
\begin{equation}\label{g31}
\frac{d}{dt}{L}(\xi,\,t) + c_1 \,f(\xi)\,\widehat{E}(\xi,\,t) \leq 0
\end{equation}
and
\begin{equation}\label{g32}
c_{3}\,\widehat{E}(\xi,\,t) \leq {L}(\xi,\,t) \leq c_{4}\,\widehat{E}(\xi,\,t),
\end{equation}
where $c_{3} =\lambda - 3c_{2} >0$ and $c_{4} =\lambda + 3c_{2} >0$. Consequently, a combination of \eqref{g31} and the second inequality in \eqref{g32} leads to, for $c=\frac{c_1}{c_4}$,  
\begin{equation}\label{g32+}
\frac{d}{dt}{L} (\xi,\,t) + c\,f(\xi)\,{L} (\xi,\,t) \leq 0.
\end{equation}
Finally, by integration \eqref{g32+} with respect to time $t$ and using \eqref{equiv} and \eqref{g32}, \eqref{11} follows with 
${\widetilde c}=\frac{c_4 \alpha_2}{c_3 \alpha_1}$.
\end{proof}
\vskip0,1truecm
\begin{theorem}\label{theorem41}
Assume that $k_2 \ne k_3$. Let $N,\,\ell\in \mathbb{N}^*$ such that $\ell\leq N$, 
\begin{equation*}
U_{0}\in H^{N}(\mathbb{R})\cap L^{1}(\mathbb{R})
\end{equation*}
and $U$ be the solution of \eqref{e11}. Then, for any $j\in\{0,\,\ldots,\,N-\ell\}$, there exists $c_0 >0$ such that 
\begin{equation}\label{33}
\Vert\partial_{x}^{j}U\Vert_{L^{2}(\mathbb{R})} \leq c_0 \,(1 + t)^{-1/12 - j/6}\,\Vert U_{0} \Vert_{L^{1}(\mathbb{R})} + c_0 \,(1 + t)^{-\ell/2}\,\Vert\partial_{x}^{j+\ell}U_{0} \Vert_{L^{2}(\mathbb{R})},\quad\forall t\in \mathbb{R}_+ .
\end{equation} 
\end{theorem}
\vskip0,1truecm
\begin{proof} From \eqref{22} we have (low and high frequences)
\begin{equation}\label{g36}
f(\xi) \geq \left\{
\begin{array}{cc}
\frac{1}{5}\,\xi^{6} & \mbox{if} \quad \vert\xi \vert\leq 1, \\
\\
\frac{1}{5}\,\xi^{-2} & \mbox{if} \quad \vert\xi\vert > 1.
\end{array}
\right.
\end{equation}
Applying Plancherel's theorem and \eqref{11}, we entail
\begin{equation}\label{g37}
\Vert\partial_{x}^{j}U\Vert_{L^{2}(\mathbb{R})}^{2} = \left\Vert\widehat{\partial_{x}^{j} U}(x,\,t)\right \Vert_{L^{2}(\mathbb{R})}^{2}  = \int_{\mathbb{R}}\xi^{2\,j}\,\vert\widehat{U}(\xi,\,t)\vert^{2} d\xi 
\end{equation}
\begin{eqnarray*}
& \leq & \widetilde{c}\int_{\mathbb{R}}\xi^{2\,j}\,e^{-\,c\,f(\xi)\,t}\,\vert\widehat{U}_0 (\xi)\vert^{2} d\xi \\
& \leq & \widetilde{c}\int_{\vert\xi\vert\leq 1}\xi^{2\,j}\,e^{-\,c\,f(\xi)\,t}\,\vert\widehat{U}_0 (\xi)\vert^{2}  d\xi + \widetilde{c}\int_{\vert\xi\vert > 1}\xi^{2\,j}\,e^{-\,c\,f(\xi)\,t}\,\vert\widehat{U}_0 (\xi)\vert^{2}\, d\xi := J_{1} + J_{2}.
\end{eqnarray*}
Using \eqref{g37+} (with $\sigma=2j$, $r=\frac{c}{5}$ and $p=6$) and \eqref{g36}, 
it follows, for the low frequency region,
\begin{equation}\label{g38}
J_{1} \leq C\,\Vert\widehat{U}_{0}\Vert_{L^{\infty}(\mathbb{R})}^{2}\int_{\vert\xi\vert\leq 1}\xi^{2\,j}\,e^{-\,\frac{c}{5}\,t\,\xi^{6}}\ d\xi \leq C\left(1 + t\right)^{-\,\frac{1}{6}(1 + 2\,j)}\Vert {U}_{0}\Vert_{L^{1}(\mathbb{R})}^{2}.
\end{equation}
For the high frequency region, using \eqref{g36}, we observe that
\begin{eqnarray*}
J_{2} & \leq & C\int_{\vert\xi\vert > 1}\vert\xi\vert^{2\,j}\,e^{-\frac{c}{5}\,t\,\xi^{-2}}\,\vert\widehat{U}(\xi,\,0)\vert^{2}\, d\xi \\
& \leq & C\ \sup_{\vert\xi\vert >1}\left\{\vert\xi\vert^{-2\,\ell}\, e^{-\frac{c}{5}\,t\,\vert\xi\vert^{-2}}\right\} \int_{\mathbb{R}}\vert\xi\vert^{2\,(j + \ell)}\,\vert\widehat{U}(\xi,\,0)\vert^{2}\ d\xi,
\end{eqnarray*}
then, using \eqref{ineq} (with $\sigma_1 =2l$, $\sigma_2 =\frac{c}{5}$ and $\sigma_3 =2$),
\begin{eqnarray}\label{g38+}
J_{2} \leq C\left(1 + t\right)^{-\,\ell}\,\Vert\,\partial_{x}^{j + \ell}U_{0}\,\Vert_{L^{2} (\mathbb{R})}^{2} , 
\end{eqnarray}
and so, by combining \eqref{g37}-\eqref{g38+}, we get \eqref{33}. 
\end{proof} 
\vskip0,1truecm
We finish this subsection by proving that the condition $k_2 \ne k_3$ is necessary for the stability of \eqref{g6} in case 
\eqref{e12} with $(\tau_1 ,\tau_2 ,\tau_3 )=(1,0,0)$. 
\vskip0,1truecm
\begin{theorem}\label{theorem410}
Assume that $k_2 = k_3$. Then $\vert\widehat{U}(\xi,\,t)\vert$ doesn't converge to zero when time $t$ goes to infinity. 
\end{theorem}
\vskip0,1truecm
\begin{proof} 
We show that, for any $\xi\in \mathbb{R}$, the matrix 
\begin{equation}\label{A012}
A:=-(-\xi^2 A_2 +i\xi A_1 +A_0)
\end{equation}
has at least a pure imaginary eigenvalue; that is   
\begin{equation*}
\forall \xi\in \mathbb{R},\,\,\exists \lambda\in \mathbb{C} :\,\,Re\,(\lambda ) =0,\quad Im\,(\lambda )\ne 0\quad\hbox{and}\quad det(\lambda I -A)=0, 
\end{equation*}
where $I$ denotes the identity matrix. From \eqref{e12} with $(\tau_1 ,\tau_2 ,\tau_3 )=(1,0,0)$ and $k_2 = k_3$, we have
\begin{equation}\label{lambdaIA}
\lambda I -A = 
\begin{pmatrix}
\lambda & -i\xi & 0 & -1 & 0 & -1 & 0\\
-ik_1 \xi & \lambda & 0 & 0 & 0 & 0 & i\gamma \xi\\
0 & 0 & \lambda & -i\xi & 0 & 0 & 0 \\
k_1 & 0 & -ik_2\xi & \lambda & 0 & 0 & 0  \\
0 & 0 & 0 & 0 & \lambda & -i\xi & 0 \\
k_1 & 0 & 0 & 0 & -ik_2 \xi & \lambda & 0\\
0 & i\gamma\xi & 0 & 0 & 0 & 0 & k_4\xi^2 +\lambda
\end{pmatrix}.
\end{equation}
A direct computaion shows that 
\begin{eqnarray}\label{detlambdaIA}
det(\lambda I -A) &=& 2k_1\lambda(\lambda^2 +k_2 \xi^2 )\left[\lambda (\lambda +k_4 \xi^2)+\gamma^2 \xi^2\right] \nonumber \\
&& +(\lambda^2 +k_2 \xi^2 )^2\left[\lambda \left(\lambda(\lambda +k_4 \xi^2)+\gamma^2 \xi^2\right)+k_1\xi^2 (\lambda +k_4 \xi^2)\right]. 
\end{eqnarray}
It is clear that, if $\xi \ne 0$, then $\lambda =i{\sqrt{k_2}} \xi$ is a pure imaginary eigenvalue of $A$.
If $\xi = 0$, then $\lambda =i{\sqrt{2k_1}}$ is a pure imaginary eigenvalue of $A$. Consequently, according to \eqref{e14} (see  \cite{tesc}), the solution of \eqref{g6} doesn't converge to zero when time $t$ goes to infinity.
\end{proof} 

\subsection{Case 2: $(\tau_1 ,\tau_2 ,\tau_3)=(0,1,0)$}

We present, first, our exponential stability result for \eqref{g6}. 
\vskip0,1truecm
\begin{lemma}\label{lemma432}
Let $\widehat{U}$ be the solution \eqref{e14} of \eqref{g6}. Then there exist $c,\,\widetilde{c}>0$ such that \eqref{11} is satisfied with  
\begin{align}\label{422}
f(\xi) = \begin{cases}
\frac{\xi^{4}}{1 +\xi^{2} +\xi^{4}}\quad &\hbox{if}\,\, k_1 =k_2 =k_3 , \\
\frac{\xi^{4}}{1 +\xi^{2} +\xi^4 +\xi^6 +\xi^{8}} \quad &\hbox{if not.}
\end{cases}
\end{align}
\end{lemma}
\vskip0,1truecm
\begin{proof}
Multiplying \eqref{e101}$_{4}$ and \eqref{e101}$_{7}$ by $i\,\xi\,\overline{\widehat{\eta}}$ and $-i\,\xi\,\overline{\widehat{y}}$, respectively, adding the resulting equations, taking the real part and using \eqref{hd2}, we get
\begin{equation}\label{equ4}
\frac{d}{dt}Re\left(i\,\xi\,\widehat{y}\, \overline{\widehat{\eta}}\right) = \gamma\xi^{2}\left(\vert\widehat{\eta}\vert^{2} - \vert\widehat{y}\vert^{2}\right) + k_4\xi^2\,Re\left(i\,\xi\,\widehat{\eta}\,\overline{\widehat{y}}\right) -k_2 \xi^2\,Re\left(\widehat{z}\, \overline{\widehat{\eta}}\right) -k_1 Re\left(i\xi\widehat{v}\, \overline{\widehat{\eta}}\right).
\end{equation}
Similarily, multiplying \eqref{e101}$_{2}$ and \eqref{e101}$_{7}$ by $\overline{\widehat{\eta}}$ and $\overline{\widehat{u}}$, respectively, adding the resulting equations, taking the real part and using \eqref{hd1}, we find
\begin{equation}\label{equ5}
\frac{d}{dt}Re\left(\widehat{u}\, \overline{\widehat{\eta}}\right) = -\gamma Re\left(i\,\xi\,\widehat{y}\,\overline{\widehat{u}}\right) -k_4 \xi^2\,Re\left(\widehat{u}\, \overline{\widehat{\eta}}\right) +k_1\,Re\left(i\xi\widehat{v}\, \overline{\widehat{\eta}}\right).
\end{equation}
Also, multiplying \eqref{e101}$_{7}$ and \eqref{e101}$_{6}$ by $i\xi\overline{\widehat{\theta}}$ and $-i\xi\overline{\widehat{\eta}}$, respectively, adding the resulting equations, taking the real part and using \eqref{hd2}, we obtain
\begin{equation}\label{equ6}
\frac{d}{dt}Re\left(i\xi\widehat{\eta}\, \overline{\widehat{\theta}}\right) = \gamma\xi^2 Re\left(\widehat{y}\,\overline{\widehat{\theta}}\right) -k_4 \xi^2\,Re\left(i\xi\widehat{\eta}\, \overline{\widehat{\theta}}\right) +k_3\xi^2\,Re\left(\widehat{\phi}\, \overline{\widehat{\eta}}\right)+k_1 Re\left(i\xi\widehat{v}\, \overline{\widehat{\eta}} \right).
\end{equation}
Let us define the functionals 
\begin{equation}\label{4g150}
{F}_0 (\xi,\,t) = Re\left[i\,\xi\,\left(\lambda_1\widehat{y}\, \overline{\widehat{z}} -\lambda_2\,\widehat{u}\, \overline{\widehat{v}} + \lambda_3\, \widehat{\theta}\, \overline{\widehat{\phi}} \right)-\lambda_4\,\xi^2\,\,\widehat{\theta}\, \overline{\widehat{v}} +\xi^2\,\widehat{y}\, \overline{\widehat{v}} \right], 
\end{equation} 
\begin{equation}\label{4g151}
{F}_1 (\xi,\,t) = \left(\frac{k_2}{k_1}\xi^2 +\lambda_1\right)Re\left(\widehat{u}\, \overline{\widehat{z}} \right), 
\end{equation}
\begin{equation}\label{4g152}
{F}_2 (\xi,\,t) = \frac{k_2}{k_1 k_3} \left(k_3 \lambda_4 \xi^2 -k_1\lambda_3\right)Re\left(i\,\xi\,\widehat{z}\, \overline{\widehat{\theta}}-\widehat{u}\, \overline{\widehat{z}} -i\,\frac{k_3}{k_2}\,\xi\,\widehat{\phi}\, \overline{\widehat{y}}\right)
\end{equation} 
and
\begin{equation}\label{4g153}
{F}_3 (\xi,\,t) = -\frac{k_2}{k_3} \left(\lambda_4 \xi^2 +\lambda_2\right)Re\left(i\,\xi\,\widehat{z}\, \overline{\widehat{\theta}}-\widehat{u}\, \overline{\widehat{z}}
-i\,\frac{k_3}{k_2}\,\xi\,\widehat{\phi}\, \overline{\widehat{y}} +\frac{k_3}{k_2} \,\widehat{u}\, \overline{\widehat{\phi}} \right),
\end{equation} 
where $\lambda_{1} ,\,\lambda_2 ,\,\lambda_{3}$ and $\lambda_{4}$ are positive constants to be fixed later.
Multiplying \eqref{eq31}-\eqref{eq35} by $\lambda_1$, $-\lambda_2$, $\lambda_3$, $\lambda_4$ and $-1$, respectively, and adding the obtained equations, it follows that
\begin{equation}\label{eq40}
\frac{d}{dt} {F}_0 (\xi,\,t) =  (\lambda_1 +1)\xi^2 \vert\widehat{y}\vert^{2} +Re\left((1-\lambda_4)\xi^2
\widehat{\theta}\, \overline{\widehat{y}}+ (\xi^2 -\lambda_2)i\xi \widehat{u}\, \overline{\widehat{y}}
-i\gamma\xi^3 \widehat{\eta}\overline{\widehat{v}}+\gamma\lambda_1\xi^2 \widehat{\eta}\overline{\widehat{z}}\right) 
\end{equation}
\begin{eqnarray*}
& & -\xi^2 \,\left((k_1 -k_1 \lambda_2 -k_1 \lambda_{4} )\vert\widehat{v}\vert^{2}
+k_2 \lambda_1 \vert\widehat{z}\vert^{2} +(\lambda_4 -\lambda_{3} )\vert\widehat{\theta}\vert^{2} + \lambda_2  \vert\widehat{u}\vert^{2} +k_3 \lambda_3\vert\widehat{\phi}\vert^{2}\right) \nonumber \\
& & +\left(k_2 \xi^2 +k_1 \lambda_1\right) \,Re\left(i\,\xi\,\widehat{z}\, \overline{\widehat{v}}\right)+\left(k_3 \lambda_4\xi^2 -k_1 \lambda_3\right) \,Re\left(i\,\xi\,\widehat{v}\, \overline{\widehat{\phi}}\right)+\left(\lambda_4\xi^2 +\lambda_2 \right) \,Re\left(i\,\xi\,\widehat{\theta}\, \overline{\widehat{u}}\right).  
\end{eqnarray*}
Multiplying \eqref{eq310} by $-\left(\frac{k_2}{k_1}\xi^2 +\lambda_1\right)$, we entail
\begin{eqnarray}\label{eq41}
\frac{d}{dt} {F}_1 (\xi,\,t) = -\left(k_2\xi^2 +k_1 \lambda_1\right)\,Re\left(i\,\xi\,\widehat{z}\, \overline{\widehat{v}}\right)
+\left(\frac{k_2}{k_1}\xi^2 +\lambda_1\right)Re\left(i\,\xi\,\widehat{y}\, \overline{\widehat{u}}\right).  
\end{eqnarray}
Multiplying \eqref{eq37} by $-\frac{k_3}{k_2}$, adding \eqref{eq36} and \eqref{eq310}, and multiplying the obtained equation by
\begin{equation*}
\frac{k_2}{k_1 k_3} \left(k_3 \lambda_4 \xi^2 -k_1\lambda_3\right), 
\end{equation*} 
we arrive at
\begin{eqnarray}\label{eq42}
\frac{d}{dt}{F}_2 (\xi,\,t) &= & \frac{k_2}{k_1 k_3} \left(k_3 \lambda_4 \xi^2 -k_1\lambda_3\right)Re\left[\left(\frac{k_3}{k_2}-1\right)\xi^2 \widehat{y}\, \overline{\widehat{\theta}}
-i\xi\widehat{y}\, \overline{\widehat{u}}+\frac{\gamma k_3}{k_2}\xi^2 \widehat{\eta}\, \overline{\widehat{\phi}}\right] 
\nonumber \\ 
& & -\left(k_3\lambda_4\xi^2 -k_1 \lambda_3\right)\,Re\left(i\,\xi\,\widehat{v}\, \overline{\widehat{\phi}}\right).
\end{eqnarray}
Similarily, adding \eqref{eq37} and \eqref{eq312}, multiplying by $-\frac{k_3}{k_2}$, adding \eqref{eq36} and \eqref{eq310}, and 
multiplying the obtained formula by $-\frac{k_2}{k_3} \left(\lambda_4 \xi^2 +\lambda_2\right)$, we deduce that
\begin{eqnarray}\label{eq43}
\frac{d}{dt} {F}_3 (\xi,\,t) & = & \frac{k_2}{k_3} \left(\lambda_4 \xi^2 +\lambda_2\right)Re\left[\left(1-\frac{k_3}{k_2}\right)\xi^2 \widehat{y}\, \overline{\widehat{\theta}}
+i\xi\widehat{y}\, \overline{\widehat{u}}-\frac{\gamma k_3}{k_2}\xi^2 \widehat{\eta}\, \overline{\widehat{\phi}}\right] \nonumber \\ 
& & -\left(\lambda_4\xi^2 + \lambda_2\right)\,Re\left(i\,\xi\,\widehat{\theta}\, \overline{\widehat{u}}\right).
\end{eqnarray}
Now, let us introduce the functional ${F}_4$
\begin{equation}\label{4g15}
{F}_4 (\xi,\,t) = {F}_0 (\xi,\,t) + {F}_1 (\xi,\,t) + {F}_2 (\xi,\,t) + {F}_3 (\xi,\,t).
\end{equation}
By combining \eqref{eq40}-\eqref{eq43}, we deduce that
\begin{eqnarray}\label{4g26+}
\frac{d}{dt}{F}_4 (\xi,\,t) &=& -\xi^{2} \,\left((k_1 -k_1 \lambda_2 -k_1 \lambda_{4} )\vert\widehat{v}\vert^{2} +k_2 \lambda_1\vert\widehat{z}\vert^{2} +(\lambda_4 -\lambda_{3} )\vert\widehat{\theta}\vert^{2} +\lambda_2 \vert\widehat{u}\vert^{2} +k_3 \lambda_3 \vert\widehat{\phi}\vert^{2}\right) \nonumber \\  
&& +F_5 (\xi,\,t),
\end{eqnarray} 
where
\begin{equation}\label{F41}
F_5 (\xi,\,t)= Re\left(I_1 \xi^2\widehat{y}\, \overline{\widehat{\theta}}+iI_2 \xi \widehat{y}\, \overline{\widehat{u}}
-I_3 \xi^2 \widehat{\eta}\, \overline{\widehat{\phi}}- i\gamma\xi^3 \widehat{\eta}\, \overline{\widehat{v}}
+\gamma\lambda_1\xi^2 \widehat{\eta}\, \overline{\widehat{z}}\right)+(\lambda_1 +1)\xi^{2}\vert\widehat{y}\vert^{2},
\end{equation} 
\begin{equation*}
I_1 =1-\lambda_4 +\frac{k_2 }{k_1 k_3}\left(\frac{k_3}{k_2} -1\right)(k_3 \lambda_4 \xi^2 -k_1 \lambda_3 )+ \left(\frac{k_2}{k_3} -1\right)\lambda_4\xi^2 +\left(\frac{k_2}{k_3} -1\right)\lambda_2 ,
\end{equation*}
\begin{equation*}
I_2 =\left[\left(\frac{k_2 }{k_3} -\frac{k_2 }{k_1}\right)\lambda_4 +\frac{k_2 }{k_1} -1\right]\xi^2 + \lambda_1 +
\left(\frac{k_2}{k_3}+1\right)\lambda_2 +\frac{k_2}{k_3}\lambda_3
\end{equation*}  
and
\begin{equation*}
I_3 =\gamma\left(1 -\frac{k_3 }{k_1}\right)\lambda_4 \xi^2 + \gamma (\lambda_2 +\lambda_3).
\end{equation*} 
Let $\lambda$ and $\lambda_5$ be positive constants, and $F$ and $L$ be the functionals  
\begin{equation}\label{F4g15}
{F}(\xi,\,t) = F_4 (\xi,\,t) +\lambda_5 Re\left(i\xi\widehat{y}\, \overline{\widehat{\eta}}\right)-\frac{1}{\gamma}I_1 
Re\left(i\xi\widehat{\eta}\, \overline{\widehat{\theta}}\right) +\frac{1}{\gamma}I_2 
Re\left(\widehat{u}\, \overline{\widehat{\eta}}\right)
\end{equation}
and 
\begin{equation}\label{4g27}
{L}(\xi,\,t) = \lambda\,\widehat{E}(\xi,\,t) +\frac{\xi^2 }{{\tilde f}(\xi) }F(\xi,\,t),
\end{equation}
where 
\begin{align}\label{4g270}
{\tilde f}(\xi) = \begin{cases}
1 +\xi^{2} +\xi^4 \quad &\hbox{if}\,\, k_1 =k_2 =k_3 , \\
1 +\xi^{2} +\xi^4 +\xi^6 +\xi^8 \quad &\hbox{if not.}
\end{cases} 
\end{align}
Multiplying \eqref{equ4}, \eqref{equ5} and \eqref{equ6} by $\lambda_5$, $\frac{1}{\gamma}I_2$ and $-\frac{1}{\gamma}I_1$, respectively, adding the obtained equations and adding \eqref{4g26+}, it appears that
\begin{eqnarray}\label{F4g26+}
\frac{d}{dt}{F} (\xi,\,t) &=& -\xi^{2} \,\left((k_1 -k_1 \lambda_2 -k_1 \lambda_{4} )\vert\widehat{v}\vert^{2} +k_2 \lambda_1\vert\widehat{z}\vert^{2} +(\lambda_4 -\lambda_{3} )\vert\widehat{\theta}\vert^{2} +\lambda_2 \vert\widehat{u}\vert^{2} +k_3 \lambda_3 \vert\widehat{\phi}\vert^{2}\right) \nonumber \\  
&& -\xi^2 \left(\gamma\lambda_5 -(\lambda_1 + 1)\right)\vert\widehat{y}\vert^{2} +F_6 (\xi,\,t),
\end{eqnarray}
where 
\begin{eqnarray}\label{7F6N1}
F_6 (\xi,\,t) &=& \gamma \lambda_5 \xi^2 \vert\widehat{\eta}\vert^{2}
+\xi^{2} Re\,\left(\gamma\lambda_1 \widehat{\eta}\, \overline{\widehat{z}} -i\gamma\xi \widehat{\eta}\, \overline{\widehat{v}} 
-I_3 \widehat{\eta}\, \overline{\widehat{\phi}} +ik_4\lambda_5 \xi \widehat{\eta}\, \overline{\widehat{y}}-k_2 \lambda_5\widehat{\eta} \, \overline{\widehat{z}}\right)-k_1 \lambda_5 Re\,\left(i\xi\widehat{v}\,\overline{\widehat{\eta}} \right)
\nonumber \\ 
&& -\frac{1}{\gamma}I_1 Re\,\left(k_3\xi^2 \widehat{\phi}\, \overline{\widehat{\eta}}+ik_1 \xi\widehat{v}\, \overline{\widehat{\eta}}-ik_4\xi^3 {\widehat{\eta}}\overline{\widehat{\theta}}\right)
+\frac{1}{\gamma}I_2 Re\,\left(ik_1 \xi\widehat{v}\, \overline{\widehat{\eta}}-k_4\xi^2 {\widehat{\eta}}\overline{\widehat{u}}\right).
\end{eqnarray}   
Noticing that, if $k_1 =k_2 =k_3$, then $I_1$, $I_2$ and $I_3$ are constants. Otherwise,
$I_1$, $I_2$ and $I_3$ are of the form $const\,\xi^2 +const$. Then, by applying Young's inequality, we see that, for any $\varepsilon_0 >0$, we have
\begin{equation}\label{F42}
F_6 (\xi,\,t) \leq \varepsilon_0 \xi^{2} \left(\vert\widehat{y}\vert^{2}+\vert\widehat{\theta}\vert^{2} +\vert\widehat{u}\vert^{2} +\vert\widehat{\phi}\vert^{2} +\vert\widehat{v}\vert^{2} +\vert\widehat{z}\vert^{2}\right) +C_{\epsilon_0 ,\lambda_1 ,\cdots,\lambda_5}{\tilde f}(\xi) \vert\widehat{\eta}\vert^{2} .
\end{equation}  
Therefore, we conclude from \eqref{F4g26+} and \eqref{F42} that 
\begin{equation} \label{F43}
\frac{d}{dt}{F}(\xi,\,t) \leq C_{\varepsilon_0 ,\lambda_1 ,\cdots,\lambda_5} {\tilde f}(\xi) \vert\widehat{\eta}\vert^{2} -\xi^2 \left(\gamma\lambda_5 -\lambda_1 - 1-\varepsilon_0\right)\vert\widehat{y}\vert^{2}
\end{equation} 
\begin{equation*}
-\xi^{2} \,\left((k_1 -k_1 \lambda_2 -k_1 \lambda_{4} -\varepsilon_0 )\vert\widehat{v}\vert^{2} +(k_2 \lambda_1 -\varepsilon_0 )\vert\widehat{z}\vert^{2} +(\lambda_4 -\lambda_{3} -\varepsilon_0)\vert\widehat{\theta}\vert^{2} +(\lambda_2 -\varepsilon_0)\vert\widehat{u}\vert^{2} +(k_3 \lambda_3 -\varepsilon_0)\vert\widehat{\phi}\vert^{2}\right) .
\end{equation*} 
We choose $0<\lambda_{1}$, $0<\lambda_{2}<1$, $0<\lambda_{3}<\lambda_{4} <1 -\lambda_2$,  $\lambda_{5} >\frac{1}{\gamma}(\lambda_{1} +1)$ and 
\begin{equation*}
0<\varepsilon_{0} <\min\left\{k_1 -k_1 \lambda_2 -k_1 \lambda_4 , k_2 \lambda_1,\lambda_4 -\lambda_3 ,
\lambda_2 ,k_3\lambda_3, \gamma\lambda_5 -\lambda_1 -1\right\}. 
\end{equation*}
Thus, using the definition \eqref{E21} of $\widehat{E}$, \eqref{F43} leads to, for some positive constant $c_1$,
\begin{equation}\label{4g27+} 
\frac{d}{dt}{F}(\xi,\,t) \leq -c_1 \xi^{2}\widehat{E} (\xi,\,t)+ C{\tilde f}(\xi)\vert\widehat{\eta}\vert^{2} .
\end{equation}
Then, from \eqref{Ep21}, \eqref{4g27} and \eqref{4g27+}, we infer that
\begin{eqnarray}\label{4g30}
\frac{d}{dt}{L}(\xi,\,t) \leq -c_1 f(\xi)\widehat{E}(\xi,\,t)-\left(k_4\,\lambda - C\right)\xi^2\vert\widehat{\eta}\vert^{2} , 
\end{eqnarray}
where $f$ is defined in \eqref{422}. On the other hand, the definitions \eqref{E21}, \eqref{F4g15} and \eqref{4g27} of
$\widehat{E}$, $F$ and $L$, respectively, imply that there exists $c_2>0$ (not depending on $\lambda$) such that  
\begin{equation*}
\left\vert {L}(\xi,\,t)-\lambda\widehat{E}(\xi,\,t)\right\vert \leq 
c_2 \frac{\xi^{2} +\vert\xi\vert^3 +\xi^{4} +\vert\xi\vert^5}{{\tilde f}(\xi)}\widehat{E}(\xi,\,t) \leq 4c_2 \widehat{E}(\xi,\,t). 
\end{equation*}
So, we choose $\lambda >\max\left\{\frac{C}{k_4},4c_2\right\}$, we get \eqref{g31} and \eqref{g32} with $c_3 =\lambda -4c_2 >0$ and
$c_4 =\lambda +4c_2 >0$. The proof can be ended as for Lemma \ref{lemma32}.
\end{proof}
\vskip0,1truecm
\begin{theorem}\label{theorem441}
Let $N,\,\ell\in \mathbb{N}^*$ such that $\ell\leq N$, 
\begin{equation*}
U_{0}\in H^{N}(\mathbb{R})\cap L^{1}(\mathbb{R})
\end{equation*}
and $U$ be the solution of \eqref{e11}. Then, for any $j\in\{0,\,\ldots,\,N-\ell\}$, there exist $c_0 ,\,{\tilde c}_0 >0$ such that, for any $t\in \mathbb{R}_+$, 
\begin{equation}\label{33case20}
\Vert\partial_{x}^{j}U\Vert_{L^{2}(\mathbb{R})} \leq c_0 \,(1 + t)^{-1/8 - j/4}\,\Vert U_{0} \Vert_{L^{1}(\mathbb{R})} + c_0 e^{-{\tilde c}_0 t}\,\Vert\partial_{x}^{j+\ell}U_{0} \Vert_{L^{2}(\mathbb{R})} \quad\hbox{if}\,\,k_1 =k_2 =k_3 , 
\end{equation} 
and 
\begin{equation}\label{33case2}
\Vert\partial_{x}^{j}U\Vert_{L^{2}(\mathbb{R})} \leq c_0 \,(1 + t)^{-1/8 - j/4}\,\Vert U_{0} \Vert_{L^{1}(\mathbb{R})} + c_0 \,(1 + t)^{-\ell/4}\,\Vert\partial_{x}^{j+\ell}U_{0} \Vert_{L^{2}(\mathbb{R})} \quad\hbox{if not}. 
\end{equation} 
\end{theorem}
\vskip0,1truecm
\begin{proof} From \eqref{422} we have (low and high frequences)
\begin{equation}
\label{4g3600}f(\xi) \geq \left\{
\begin{array}{cc}
\frac{1}{3}\,\xi^{4} & \mbox{if} \quad \vert\xi \vert\leq 1, \\
\\
\frac{1}{3} & \mbox{if} \quad \vert\xi\vert > 1
\end{array}
\right. \quad\hbox{if}\,\,k_1 =k_2 =k_3 , 
\end{equation}
and
\begin{equation}
\label{4g36}f(\xi) \geq \left\{
\begin{array}{cc}
\frac{1}{5}\,\xi^{4} & \mbox{if} \quad \vert\xi \vert\leq 1, \\
\\
\frac{1}{5}\,\xi^{-4} & \mbox{if} \quad \vert\xi\vert > 1
\end{array}
\right. \quad\hbox{if not}.
\end{equation}
The proof of \eqref{33case2} is identical to the one of Theorem \ref{theorem41} by using \eqref{4g36} and applying \eqref{g37+} (with $\sigma =2j$, $r=\frac{c}{5}$ and $p=4$) and \eqref{ineq} (with $\sigma_1 =2l$, $\sigma_2 =\frac{c}{5}$ and $\sigma_3 =4$). To get \eqref{33case20}, noticing that the low frequencies can be treated as for \eqref{33case2}. For the high frequencies, we have just to remark that \eqref{4g3600} implies that
\begin{eqnarray*}
\int_{\vert\xi\vert > 1}\vert\xi\vert^{2\,j}\,e^{-cf(\xi)t}\,\vert\widehat{U}(\xi,\,0)\vert^{2}\ d\xi& \leq & \int_{\vert\xi\vert > 1}\vert\xi\vert^{2\,j}\,e^{-\frac{c}{3}\,t}\,\vert\widehat{U}(\xi,\,0)\vert^{2}\ d\xi  \nonumber \\
& \leq & \sup_{\vert\xi\vert >1}\left\{\vert\xi\vert^{-2\,\ell}\, e^{-\frac{c}{3}\,t}\right\} \int_{\mathbb{R}}\vert\xi\vert^{2\,(j + \ell)}\,\vert\widehat{U}(\xi,\,0)\vert^{2}\ d\xi  \nonumber \\
& \leq & e^{-\frac{c}{3}\,t}\,\Vert\,\partial_{x}^{j + \ell}U_{0}\,\Vert_{L^{2} (\mathbb{R})}^{2} ,
\end{eqnarray*}
so \eqref{33case20} holds true with ${\tilde c}_0 =\frac{c}{6}$.
\end{proof}

\subsection{Case 3: $(\tau_1 ,\tau_2 ,\tau_3)=(0,0,1)$}

In this case, we prove the same stability results for \eqref{g6} and \eqref{e11} that given in Subsection 3.2, and moreover, the proofs are very similar. 
\vskip0,1truecm
\begin{lemma}\label{lemma532}
The result of Lemma \ref{lemma432} holds true also when $(\tau_1 ,\tau_2 ,\tau_3)=(0,0,1)$. 
\end{lemma}
\vskip0,1truecm
\begin{proof}
Multiplying \eqref{e101}$_{6}$ and \eqref{e101}$_{7}$ by $i\,\xi\,\overline{\widehat{\eta}}$ and $-i\,\xi\,\overline{\widehat{\theta}}$, respectively, adding the resulting equations, taking the real part and using \eqref{hd2}, we get
\begin{equation}\label{equ7}
\frac{d}{dt}Re\left(i\,\xi\,\widehat{\theta}\, \overline{\widehat{\eta}}\right) = \gamma\xi^{2}\left(\vert\widehat{\eta}\vert^{2} - \vert\widehat{\theta}\vert^{2}\right) + k_4\xi^2\,Re\left(i\,\xi\,\widehat{\eta}\,\overline{\widehat{\theta}}\right) -k_3 \xi^2\,Re\left(\widehat{\phi}\, \overline{\widehat{\eta}}\right) -k_1 Re\left(i\xi\widehat{v}\, \overline{\widehat{\eta}}\right).
\end{equation}
Similarily, multiplying \eqref{e101}$_{2}$ and \eqref{e101}$_{7}$ by $\overline{\widehat{\eta}}$ and $\overline{\widehat{u}}$, respectively, adding the resulting equations, taking the real part and using \eqref{hd1}, we find
\begin{equation}\label{equ8}
\frac{d}{dt}Re\left(\widehat{u}\, \overline{\widehat{\eta}}\right) = -\gamma Re\left(i\,\xi\,\widehat{\theta}\,\overline{\widehat{u}}\right) -k_4 \xi^2\,Re\left(\widehat{u}\, \overline{\widehat{\eta}}\right) +k_1\,Re\left(i\xi\widehat{v}\, \overline{\widehat{\eta}}\right).
\end{equation}
Also, multiplying \eqref{e101}$_{7}$ and \eqref{e101}$_{4}$ by $i\xi\overline{\widehat{y}}$ and $-i\xi\overline{\widehat{\eta}}$, respectively, adding the resulting equations, taking the real part and using \eqref{hd2}, we obtain
\begin{equation}\label{equ9}
\frac{d}{dt}Re\left(i\xi\widehat{\eta}\, \overline{\widehat{y}}\right) = \gamma\xi^2 Re\left(\widehat{\theta}\,\overline{\widehat{y}}\right) -k_4 \xi^2\,Re\left(i\xi\widehat{\eta}\, \overline{\widehat{y}}\right) +k_2\xi^2\,Re\left(\widehat{z}\, \overline{\widehat{\eta}}\right)+k_1 Re\left(i\xi\widehat{v}\, \overline{\widehat{\eta}} \right).
\end{equation}
After, we define the functionals 
\begin{equation}\label{5g150}
{F}_0 (\xi,\,t) = Re\left[i\,\xi\,\left(\lambda_1\widehat{y}\, \overline{\widehat{z}} -\lambda_2\,\widehat{u}\, \overline{\widehat{v}} + \lambda_3\,\widehat{\theta}\, \overline{\widehat{\phi}}\right)+\lambda_4\,\xi^2\,\,\widehat{\theta}\, \overline{\widehat{v}} -\xi^2\,\widehat{y}\, \overline{\widehat{v}} \right], 
\end{equation} 
\begin{equation}\label{5g151}
{F}_1 (\xi,\,t) = \left(\frac{k_3}{k_1}\lambda_4 \xi^2 +\lambda_3\right)Re\left(\widehat{u}\, \overline{\widehat{\phi}} \right), 
\end{equation}
\begin{equation}\label{5g152}
{F}_2 (\xi,\,t) = -\frac{1}{k_1} \left(k_2 \xi^2 -k_1\lambda_1\right)Re\left(i\,\xi\,\widehat{z}\, \overline{\widehat{\theta}}-i\frac{k_3}{k_2}\xi\widehat{\phi}\, \overline{\widehat{y}} +\frac{k_3}{k_2}\widehat{u}\, \overline{\widehat{\phi}}\right)
\end{equation} 
and
\begin{equation}\label{5g153}
{F}_3 (\xi,\,t) = \left(\xi^2 +\lambda_2\right)Re\left(i\,\xi\,\widehat{z}\, \overline{\widehat{\theta}}-i\frac{k_3}{k_2}\xi\widehat{\phi}\, \overline{\widehat{y}} +\frac{k_3}{k_2}\widehat{u}\, \overline{\widehat{\phi}} - \widehat{u}\, \overline{\widehat{z}}\right),
\end{equation} 
where $\lambda_{1} ,\,\lambda_2 ,\,\lambda_{3}$ and $\lambda_{4}$ are positive constants to be fixed later.
Multiplying \eqref{eq31}-\eqref{eq34} by $\lambda_1$, $-\lambda_2$, $\lambda_3$ and $-\lambda_4$, respectively, adding the obtained equations and adding \eqref{eq35}, we infer that
\begin{equation}\label{eq50}
\frac{d}{dt} {F}_0 (\xi,\,t) =Re\left[i\xi\left(\left(\lambda_2 -\lambda_4 \xi^2 \right)\widehat{\theta}\, \overline{\widehat{u}} -i\xi\gamma\lambda_3 \widehat{\eta}\, \overline{\widehat{\phi}}\right)+(\lambda_4 -1)\xi^2 \widehat{\theta}\, \overline{\widehat{y}}-i\gamma \lambda_4 \xi^3 \widehat{\eta}\, \overline{\widehat{v}}\right]  
\end{equation}
\begin{eqnarray*}
& & + (\lambda_3 +\lambda_4)\xi^2\vert\widehat{\theta}\vert^{2} -\xi^2 \,\left((k_1 \lambda_4 -k_1\lambda_{2} -k_1 )\vert\widehat{v}\vert^{2}
+k_2 \lambda_1\vert\widehat{z}\vert^{2} +(1-\lambda_1 )\vert\widehat{y}\vert^{2} +\lambda_2 \vert\widehat{u}\vert^{2} +k_3 \lambda_3 \vert\widehat{\phi}\vert^{2}\right) \nonumber \\
& & +(k_2 \xi^2 -k_1 \lambda_1) \,Re\left(i\,\xi\,\widehat{v}\, \overline{\widehat{z}}\right)+(k_3 \lambda_4\xi^2 +k_1 \lambda_3) \,Re\left(i\,\xi\,\widehat{\phi}\, \overline{\widehat{v}}\right)+(\xi^2 +\lambda_2 ) \,Re\left(i\,\xi\,\widehat{y}\, \overline{\widehat{u}}\right).  
\end{eqnarray*}
Multiplying \eqref{eq312} by $-\left(\frac{k_3}{k_1}\lambda_4 \xi^2 +\lambda_3\right)$, we arrive at
\begin{eqnarray}\label{eq51}
\frac{d}{dt} {F}_1 (\xi,\,t) = \left(\frac{k_3}{k_2}\lambda_4\xi^2 +\lambda_3\right)\,Re\left(i\,\xi\,\widehat{\theta}\, \overline{\widehat{u}}\right) -\left(k_3\lambda_4\xi^2 +k_1 \lambda_3\right)\,Re\left(i\,\xi\,\widehat{\phi}\, \overline{\widehat{v}}\right).  
\end{eqnarray}
Adding \eqref{eq37} and \eqref{eq312}, multiplying the obtained equation by $-\frac{k_3}{k_2}$, adding \eqref{eq36} and multiplying the reuslting equation by $-\left(\frac{k_2}{k_1} \xi^2 -\lambda_1\right)$, it follows that
\begin{eqnarray}\label{eq52}
\frac{d}{dt}{F}_2 (\xi,\,t) & = & \left(\frac{k_2}{k_1}\xi^2 -\lambda_1\right)Re\left[i\xi\left(-\frac{k_3}{k_2}
\widehat{\theta}\, \overline{\widehat{u}} -i\gamma\xi\widehat{\eta}\, \overline{\widehat{z}}\right)+\left(1-\frac{k_3}{k_2}\right)\xi^2 \widehat{\theta}\, \overline{\widehat{y}}\right]  \nonumber \\ 
& & -\left(k_2\xi^2 -k_1 \lambda_1\right)\,Re\left(i\,\xi\,\widehat{v}\, \overline{\widehat{z}}\right).
\end{eqnarray}
Similarily, adding \eqref{eq37} and \eqref{eq312}, multiplying the obtained equation by $-\frac{k_3}{k_2}$, adding \eqref{eq36}
and \eqref{eq310}, and multiplying the reuslting equation by $\xi^2 +\lambda_2$, we entail
\begin{eqnarray}\label{eq53}
\frac{d}{dt} {F}_3 (\xi,\,t) & = & \left(\xi^2 +\lambda_2\right)Re\left[i\xi\left(\frac{k_3}{k_2}
\widehat{\theta}\, \overline{\widehat{u}} +i\gamma\xi\widehat{\eta}\, \overline{\widehat{z}}\right)+\left(\frac{k_3}{k_2} -1\right) \xi^2 \widehat{\theta}\, \overline{\widehat{y}}\right]\nonumber \\ 
& & -\left(\xi^2 + \lambda_2\right)\,Re\left(i\,\xi\,\widehat{y}\, \overline{\widehat{u}}\right).
\end{eqnarray}
Let $F_4$ the functional defined by \eqref{4g15}. A combination of \eqref{eq50}-\eqref{eq53} implies that 
\begin{eqnarray} \label{5g26+}
\frac{d}{dt}{F}_4 (\xi,\,t) &=& -\xi^2 \,\left((k_1 \lambda_4 -k_1\lambda_{2} -k_1 )\vert\widehat{v}\vert^{2}+k_2 \lambda_1\vert\widehat{z}\vert^{2} +(1-\lambda_1 )\vert\widehat{y}\vert^{2} + \lambda_2 \vert\widehat{u}\vert^{2} +k_3 \lambda_3 \vert\widehat{\phi}\vert^{2}\right) \nonumber \\ 
&& +F_5 (\xi,\,t),
\end{eqnarray} 
where
\begin{equation}\label{F410}
F_5 (\xi,\,t)= Re\left(iI_1 \xi\widehat{\theta}\, \overline{\widehat{u}}-I_2 \xi^2 \widehat{\eta}\, \overline{\widehat{z}}
+I_3 \xi^2\widehat{\theta}\, \overline{\widehat{y}}-i\gamma \lambda_4 \xi^3 \widehat{\eta}\, \overline{\widehat{v}}
+\gamma\lambda_3\xi^2 \widehat{\eta}\, \overline{\widehat{\phi}}\right)+(\lambda_3 +\lambda_4 )\xi^{2}\vert\widehat{\theta}\vert^{2},
\end{equation} 
\begin{equation*}
I_1 =\left[\left(\frac{k_3 }{k_2} -1\right)\lambda_4 +\frac{k_3 }{k_2}-\frac{k_3 }{k_1}\right]\xi^2 +\frac{k_3}{k_2} \lambda_1 +
\left(\frac{k_3 }{k_2} +1\right)\lambda_2 +\lambda_3 ,
\end{equation*}
\begin{equation*}
I_2 =\gamma \left[\left(1-\frac{k_2 }{k_1}\right)\xi^2 +\lambda_1 +\lambda_2 \right]
\end{equation*}  
and
\begin{equation*}
I_3 =\left(\frac{k_3 }{k_2} -1\right)(\xi^2 +\lambda_2) +\left(1-\frac{k_3 }{k_2} \right)\left(\frac{k_2 }{k_1}\xi^2 -\lambda_1\right)+\lambda_4 -1.
\end{equation*} 
Let $\lambda$ and $\lambda_5$ be positive constants and $L$ be the functional defined by \eqref{4g27}, where ${\tilde f}$ is defined by \eqref{4g270} and $F$ is given by 
\begin{equation}\label{FF4g15}
{F}(\xi,\,t) = F_4 (\xi,\,t) +\lambda_5 Re\left(i\xi\widehat{\theta}\, \overline{\widehat{\eta}}\right)+\frac{1}{\gamma}I_1 
Re\left(\widehat{u}\, \overline{\widehat{\eta}}\right) -\frac{1}{\gamma}I_3 
Re\left(i\xi\widehat{\eta}\, \overline{\widehat{y}}\right).
\end{equation}
Multiplying \eqref{equ7}-\eqref{equ9} by $\lambda_5$, $\frac{1}{\gamma}I_1$ and $-\frac{1}{\gamma}I_3$, respectively, adding the obtained equations and adding \eqref{5g26+}, we find
\begin{eqnarray}\label{FF4g26+}
\frac{d}{dt}{F} (\xi,\,t) &=& -\xi^{2} \,\left((k_1 \lambda_4 -k_1\lambda_{2} -k_1 )\vert\widehat{v}\vert^{2}+k_2 \lambda_1\vert\widehat{z}\vert^{2} +(1-\lambda_1 )\vert\widehat{y}\vert^{2} + \lambda_2 \vert\widehat{u}\vert^{2} +k_3 \lambda_3 \vert\widehat{\phi}\vert^{2}\right) \nonumber \\  
&& -\left(\gamma\lambda_5 -\lambda_3 -\lambda_4\right)\xi^2 \vert\widehat{\theta}\vert^{2} +F_6 (\xi,\,t),
\end{eqnarray}
where 
\begin{eqnarray}\label{7F6N}
F_6 (\xi,\,t) &=& \gamma \lambda_5 \xi^2 \vert\widehat{\eta}\vert^{2}
-\xi^{2} Re\left(I_2 \widehat{\eta}\, \overline{\widehat{z}}+i\gamma \lambda_4 \xi \widehat{\eta}\, \overline{\widehat{v}}
-\gamma\lambda_3\widehat{\eta}\, \overline{\widehat{\phi}}-ik_4 \lambda_5 \xi \widehat{\eta}\, \overline{\widehat{\theta}}+k_3\lambda_5\widehat{\phi}\, \overline{\widehat{\eta}}\right) -k_1 \lambda_5 Re\left(i\xi \widehat{v}\, \overline{\widehat{\eta}}\right) \nonumber \\ 
&& -\frac{1}{\gamma}I_1 Re\,\left(k_4\xi^2 \widehat{\eta}\, \overline{\widehat{u}}-ik_1 \xi\widehat{v}\, \overline{\widehat{\eta}}\right)+\frac{1}{\gamma}I_3 Re\,\left(ik_4 \xi^3\widehat{\eta}\, \overline{\widehat{y}}-k_2\xi^2 {\widehat{\eta}}\overline{\widehat{z}}-ik_1\xi {\widehat{v}}\overline{\widehat{\eta}}\right).
\end{eqnarray}   
We remark that, if $k_1 =k_2 =k_3$, then $I_1$, $I_2$ and $I_3$ are constants. Otherwise,
$I_1$, $I_2$ and $I_3$ are of the form $const\,\xi^2 +const$. Then, by applying Young's inequality, we see that, for any $\varepsilon_0 >0$, we have
\begin{equation}\label{FF42}
F_6 (\xi,\,t) \leq \varepsilon_0 \xi^{2} \left(\vert\widehat{y}\vert^{2}+\vert\widehat{\theta}\vert^{2} +\vert\widehat{u}\vert^{2} +\vert\widehat{\phi}\vert^{2} +\vert\widehat{v}\vert^{2} +\vert\widehat{z}\vert^{2}\right) +C_{\epsilon_0 ,\lambda_1 ,\cdots,\lambda_5}{\tilde f}(\xi) \vert\widehat{\eta}\vert^{2} ,
\end{equation}
where ${\tilde f}$ is defined in \eqref{4g270}. Therefore, we conclude from \eqref{FF4g26+} and \eqref{FF42} that 
\begin{equation} \label{FF43}
\frac{d}{dt}{F}(\xi,\,t) \leq C_{\varepsilon_0 ,\lambda_1 ,\cdots,\lambda_5} {\tilde f}(\xi) \vert\widehat{\eta}\vert^{2} - \left(\gamma\lambda_5 -\lambda_3 - \lambda_4 -\varepsilon_0\right)\xi^2\vert\widehat{\theta}\vert^{2}
\end{equation} 
\begin{equation*}
-\xi^{2} \,\left((k_1 \lambda_4 -k_1\lambda_{2} -k_1 -\varepsilon_0 )\vert\widehat{v}\vert^{2}+(k_2 \lambda_1 -\varepsilon_0 ) \vert\widehat{z}\vert^{2} +(1-\lambda_1 -\varepsilon_0 )\vert\widehat{y}\vert^{2} + (\lambda_2 -\varepsilon_0 ) \vert\widehat{u}\vert^{2} +(k_3 \lambda_3 -\varepsilon_0 )\vert\widehat{\phi}\vert^{2}\right) .
\end{equation*}
We choose $0<\lambda_{3}$, $0<\lambda_{1}<1$, $\lambda_{4} >1$, $0<\lambda_{2}<\lambda_{4} -1$, 
$\lambda_{5} >\frac{1}{\gamma} (\lambda_3 +\lambda_4 )$ and 
\begin{equation*}
0<\varepsilon_{0} <\min\left\{k_3\lambda_3 , \lambda_2,
1-\lambda_1 ,k_2\lambda_1 , k_1\lambda_4 -k_1\lambda_2 -k_1, \gamma \lambda_5 -\lambda_3 -\lambda_4\right\}. 
\end{equation*}
Then, using the definition \eqref{E21} of $\widehat{E}$, \eqref{FF43} imlies \eqref{4g27+}, and then \eqref{4g30} holds true. 
Consequentely, the proof can be ended as for Lemma \ref{lemma432}.
\end{proof}
\vskip0,1truecm
\begin{theorem}\label{theorem541}
The stability result given in Theorem \ref{theorem441} is satisfied when $(\tau_1 ,\tau_2 ,\tau_3)=(0,0,1)$.
\end{theorem}
\vskip0,1truecm
\begin{proof} 
The proof is identical to the one of Theorem \ref{theorem441}.
\end{proof}

\section{Stability: Cattaneo law \eqref{s22}}

This section concerns the stability of \eqref{e11} in case of Cattaneo law \eqref{s22}. We will prove \eqref{estim1}, \eqref{estim2} and \eqref{estim3}. Moreover, we prove that \eqref{g6} is not stable when $(\tau_1 ,\tau_2 ,\tau_3)=(1,0,0)$ and 
$k_2 =k_3$.   
\vskip0,1truecm
First, observe that \eqref{e101}$_{1}$-\eqref{e101}$_{6}$ are identical to \eqref{e102}$_{1}$-\eqref{e102}$_{6}$, and
\eqref{e101}$_{7}$ with $k_4\xi^2\widehat{\eta}$ replaced by $ik_4 \xi \widehat{q}$ is equal to \eqref{e102}$_{7}$. So \eqref{eq31}-\eqref{eq312} are still valide. Moreover, \eqref{equ1}-\eqref{equ3}, \eqref{equ4}-\eqref{equ6} and \eqref{equ7}-\eqref{equ9} are satisfied with $ik_4 \xi\widehat{q}$ instead of $k_4\xi^2\widehat{\eta}$. On the other hand, we prove the next expressions, which take in consideration the last equation in \eqref{e102}. 
\vskip0,1truecm
Multiplying \eqref{e102}$_{7}$ and \eqref{e102}$_{8}$ by $i\xi\overline{\widehat{q}}$ and $-i\xi\overline{\widehat{\eta}}$, respectively, adding the resulting equations, taking the real part and using \eqref{hd2}, we find
\begin{equation}\label{7equ1}
\frac{d}{dt}Re\left(i\xi\widehat{\eta}\, \overline{\widehat{q}}\right) = k_4 \xi^2 \left(\vert\widehat{q}\vert^{2} - \vert\widehat{\eta}\vert^{2}\right) +k_5 Re\left(i\xi\widehat{q}\, \overline{\widehat{\eta}}\right) +\gamma\xi^2\,Re\left((\tau_1\widehat{u}+\tau_2\widehat{y} +\tau_3\widehat{\theta})\, \overline{\widehat{q}}\right).
\end{equation}
Multiplying \eqref{e102}$_{1}$ and \eqref{e102}$_{8}$ by $i\,\xi\,\overline{\widehat{q}}$ and $-i\,\xi\,\overline{\widehat{v}}$, respectively, adding the resulting equations, taking the real part and using \eqref{hd2}, we get
\begin{equation}\label{7equ2}
\frac{d}{dt}Re\left(i\,\xi\,\widehat{v}\, \overline{\widehat{q}}\right) = -\xi^{2} Re\left(\widehat{u}\,\overline{\widehat{q}}\right)
+Re\left(i\xi\widehat{y}\,\overline{\widehat{q}}+i\xi\widehat{\theta}\,\overline{\widehat{q}}+ik_5 \xi\widehat{q}\,\overline{\widehat{v}}\right) -k_4 \xi^2\,Re\left(\widehat{v}\, \overline{\widehat{\eta}}\right).
\end{equation}
Similarily, using the multipliers $\overline{\widehat{q}}$ and $\overline{\widehat{v}}$ instead of $i\,\xi\,\overline{\widehat{q}}$ and $-i\,\xi\,\overline{\widehat{v}}$, respectively, we obtain
\begin{equation}\label{7equ3}
\frac{d}{dt}Re\left(\widehat{v}\, \overline{\widehat{q}}\right) = Re\left(i\xi\widehat{u}\,\overline{\widehat{q}}\right)
+Re\left(\widehat{y}\,\overline{\widehat{q}}+\widehat{\theta}\,\overline{\widehat{q}}+k_5\widehat{q}\,\overline{\widehat{v}}\right) +k_4 \,Re\left(i\xi\widehat{v}\, \overline{\widehat{\eta}}\right).
\end{equation}
Multiplying \eqref{e102}$_{5}$ and \eqref{e102}$_{8}$ by $i\,\xi\,\overline{\widehat{q}}$ and $-i\,\xi\,\overline{\widehat{\phi}}$, respectively, adding the resulting equations, taking the real part and using \eqref{hd2}, we arrive at
\begin{equation}\label{7equ4}
\frac{d}{dt}Re\left(i\,\xi\,\widehat{\phi}\, \overline{\widehat{q}}\right) = -\xi^{2} Re\left(\widehat{\theta}\,\overline{\widehat{q}}\right)
+k_5 Re\left(i\xi\widehat{q}\,\overline{\widehat{\phi}}\right) -k_4 \xi^2\,Re\left(\widehat{\phi}\, \overline{\widehat{\eta}}\right).
\end{equation}
Similarily, using the multipliers $\overline{\widehat{q}}$ and $\overline{\widehat{\phi}}$ instead of $i\,\xi\,\overline{\widehat{q}}$ and $-i\,\xi\,\overline{\widehat{\phi}}$, respectively, we entail
\begin{equation}\label{7equ5}
\frac{d}{dt}Re\left(\widehat{\phi}\, \overline{\widehat{q}}\right) = Re\left(i\xi\widehat{\theta}\,\overline{\widehat{q}}\right)
+k_5 Re\left(\widehat{q}\,\overline{\widehat{\phi}}\right) +k_4 \,Re\left(i\xi\widehat{\phi}\, \overline{\widehat{\eta}}\right).
\end{equation}
Also, multiplying \eqref{e102}$_{3}$ and \eqref{e102}$_{8}$ by $i\,\xi\,\overline{\widehat{q}}$ and $-i\,\xi\,\overline{\widehat{z}}$, respectively, adding the resulting equations, taking the real part and using \eqref{hd2}, we infer that
\begin{equation}\label{7equ6}
\frac{d}{dt}Re\left(i\,\xi\,\widehat{z}\, \overline{\widehat{q}}\right) = -\xi^{2} Re\left(\widehat{y}\,\overline{\widehat{q}}\right)
+k_5 Re\left(i\xi\widehat{q}\,\overline{\widehat{z}}\right) -k_4 \xi^2\,Re\left(\widehat{z}\, \overline{\widehat{\eta}}\right).
\end{equation}
Similarily, using the multipliers $\overline{\widehat{q}}$ and $\overline{\widehat{z}}$ instead of $i\,\xi\,\overline{\widehat{q}}$ and $-i\,\xi\,\overline{\widehat{z}}$, respectively, it appears that
\begin{equation}\label{7equ7}
\frac{d}{dt}Re\left(\widehat{z}\, \overline{\widehat{q}}\right) = Re\left(i\xi\widehat{y}\,\overline{\widehat{q}}\right)
+k_5 Re\left(\widehat{q}\,\overline{\widehat{z}}\right) +k_4 \,Re\left(i\xi\widehat{z}\, \overline{\widehat{\eta}}\right).
\end{equation}

\subsection{Case 1: $(\tau_1 ,\tau_2 ,\tau_3)=(1,0,0)$}

As in Subsection 3.1, we start by presenting the exponential stability result of \eqref{g6} in the next lemma. 
\vskip0,1truecm
\begin{lemma}\label{lemma732}
The result of Lemma \ref{lemma32} is satisfied in case \eqref{e102} when $k_2 \ne k_3$ and $(\tau_1 ,\tau_2 ,\tau_3)=(1,0,0)$.
\end{lemma}
\vskip0,1truecm
\begin{proof}
We use the arguments used in Subsection 3.1. We define the functional ${F}_0$ by \eqref{F0} and we get \eqref{g26+} (because we used only the first six equations in \eqref{e101} which are the same in \eqref{e102}). We consider $F_1$ and $F$ defined by \eqref{F} and \eqref{FF1}, and we find \eqref{FFprime1} with $k_4\xi\widehat{\eta}$ replaced by $ik_4\widehat{q}$. We put, for $\lambda_6 >0$, 
\begin{eqnarray*}
{\tilde F} (\xi,\,t) = {F} (\xi,\,t)+ \lambda_6 \xi^4 Re\left(i\xi\widehat{\eta}\, \overline{\widehat{q}} \right)+\frac{1}{k_4}I_5\xi^2 
Re\left(i\xi\widehat{v}\, \overline{\widehat{q}} \right) +\frac{1}{k_4}I_6\xi^4 Re\left(\widehat{\phi}\, \overline{\widehat{q}} \right)+\frac{1}{k_4}I_7\xi^4 Re\left(\widehat{z}\, \overline{\widehat{q}} \right),
\end{eqnarray*}
where
\begin{equation*}
I_5 =(\gamma\lambda_2 -k_1\lambda_5 )\xi^2 -\frac{k_1}{\gamma}(I_3 +I_4),\quad I_6 =\frac{\gamma}{k_1} I_2 -\frac{k_3}{\gamma}I_3\quad\hbox{and}\quad I_7 =\frac{\gamma}{k_1} I_1 -\frac{k_2}{\gamma}I_4 . 
\end{equation*}
Multiplying \eqref{7equ1}, \eqref{7equ2}, \eqref{7equ5} and \eqref{7equ7} by $\lambda_6 \xi^4$, $\frac{1}{k_4}I_5\xi^2$,
$\frac{1}{k_4}I_6\xi^4$ and $\frac{1}{k_4}I_7\xi^4$, respectively, adding the obtained equations, adding \eqref{FFprime1} and applying Young's inequality for the terms depending on $\widehat{q}$, we find, for any $\varepsilon_0 >0$,
\begin{eqnarray} \label{7Fprime2}
\frac{d}{dt}{\tilde F} (\xi,\,t) &\leq & -(k_2\lambda_1 -\varepsilon_0 ) \xi^{6}\,\vert\widehat{z}\vert^{2} 
-(k_3\lambda_3 -\varepsilon_0 )\xi^{6}\,\vert\widehat{\phi}\vert^{2} -\left(1-\lambda_1 -\varepsilon_0 \right)\xi^{6}\,\vert\widehat{y}\vert^{2} \nonumber \\ 
& & -\left(\lambda_4 -\lambda_3 -\varepsilon_0 \right)\xi^{6}\,\vert\widehat{\theta}\vert^{2} -(k_1\lambda_2 -k_1\lambda_4 -k_1  -\varepsilon_0 )\xi^{6} \,\vert\widehat{v}\vert^{2} -\left(\gamma \lambda_5 -\lambda_2 -\varepsilon_0 \right)\xi^{6}\,\vert\widehat{u}\vert^{2} \nonumber \\
& & -\left(k_4 \lambda_6 -\gamma\lambda_5 -\varepsilon_0 \right)\xi^{6}\,\vert\widehat{\eta}\vert^{2}+C_{\varepsilon_0 ,\lambda_1 ,\cdots ,\lambda_6} (1+\xi^2 +\xi^4 +\xi^6 +\xi^8) \vert\widehat{q}\vert^{2} .
\end{eqnarray}
We choose $0<\lambda_{1}<1$, $\lambda_{2} >1$,  $\lambda_{5} >\frac{1}{\gamma}\lambda_{2}$, 
$\lambda_{6} >\frac{\gamma}{k_4}\lambda_{5}$, $0<\lambda_{3}<\lambda_{4}<\lambda_{2} -1$ and 
\begin{equation*}
0<\varepsilon_{0} <\min\left\{k_2\lambda_1 ,k_3\lambda_3 ,1- \lambda_1 ,\lambda_4 - \lambda_3 , k_1\lambda_2 -k_1\lambda_4 -k_1 , \gamma \lambda_5 -\lambda_2, k_4 \lambda_6 -\gamma\lambda_5\right\}. 
\end{equation*} 
Hence, using the definition \eqref{E22} of $\widehat{E}$, \eqref{7Fprime2} leads to, for some positive constant $c_1$,
\begin{eqnarray}\label{7g27+} 
\frac{d}{dt}{\tilde F} (\xi,\,t) \leq -c_1 \xi^{6}\widehat{E} (\xi,\,t)+ C \,\left(1 + \xi^2 +\xi^{4} +\xi^{6} +\xi^{8} \right)\vert\widehat{q}\vert^{2} .
\end{eqnarray}
So we consider ${L}$  given by \eqref{g27}, with ${\tilde F}$ instead of $F$, and use \eqref{Ep22} to find  
\begin{eqnarray}\label{7g30}
\frac{d}{dt}{L}(\xi,\,t) \leq -c_1 f(\xi)\widehat{E}(\xi,\,t)-\left(k_5\,\lambda - C\right)\vert\widehat{q}\vert^{2} , 
\end{eqnarray}
where $f$ is defined by \eqref{22}. Finally, the proof can be finished exactely as in the proof of \eqref{11}. 
\end{proof}
\vskip0,1truecm
\begin{theorem}\label{theorem741}
The result of Theorem \ref{theorem41} is satisfied in case \eqref{e102} when $k_2 \ne k_3$ and $(\tau_1 ,\tau_2 ,\tau_3)=(1,0,0)$.
\end{theorem}
\vskip0,1truecm
\begin{proof} The proof is identical to the one of Theorem \ref{theorem41}.
\end{proof} 
\vskip0,1truecm
The third result of this subsection says that \eqref{g6} is not stable if $k_2 = k_3$.
\vskip0,1truecm
\begin{theorem}\label{theorem7410}
Assume that $k_2 = k_3$. Then $\vert\widehat{U}(\xi,\,t)\vert$ doesn't converge to zero when time $t$ goes to infinity. 
\end{theorem}
\vskip0,1truecm
\begin{proof} 
As in Subsection 3.1, we show that, for any $\xi\in \mathbb{R}$, the matrix \eqref{A012} 
has at least a pure imaginary eigenvalue. From \eqref{e13} with $(\tau_1 ,\tau_2 ,\tau_3 )=(1,0,0)$ and $k_2 = k_3$, we have
\begin{equation}\label{lambdaIA2}
\lambda I -A = 
\begin{pmatrix}
\lambda & -i\xi & 0 & -1 & 0 & -1 & 0 & 0\\
-ik_1 \xi & \lambda & 0 & 0 & 0 & 0 & i\gamma \xi & 0\\
0 & 0 & \lambda & -i\xi & 0 & 0 & 0 & 0 \\
k_1 & 0 & -ik_2\xi & \lambda & 0 & 0 & 0 & 0 \\
0 & 0 & 0 & 0 & \lambda & -i\xi & 0 & 0 \\
k_1 & 0 & 0 & 0 & -ik_2 \xi & \lambda & 0 & 0\\
0 & i\gamma\xi & 0 & 0 & 0 & 0 & \lambda & ik_4\xi\\
0 & 0 & 0 & 0 & 0 & 0 & ik_4\xi & k_5 +\lambda
\end{pmatrix} .
\end{equation}
A direct computaion shows that 
\begin{eqnarray}\label{detlambdaIA2}
det(\lambda I -A) &=& 2k_1\lambda(\lambda +k_5 )(\lambda^2 +k_2 \xi^2 )\left(\lambda^2 +\gamma^2 \xi^2\right) +\lambda (\lambda +k_5 )
(\lambda^2 +k_2 \xi^2 )^2 \left(\lambda^2 +(k_1 +\gamma^2 )\xi^2\right)\nonumber\\
&& -ik_4\xi (\lambda^2 +k_2 \xi^2 )\left[ik_1 k_4 \xi^3 (\lambda^2 +k_2 \xi^2 )+ik_4 \lambda^2\xi (\lambda +k_1 )+
ik_2 k_4 \lambda^2\xi^3 +ik_1 k_4 \lambda\xi\right]. 
\end{eqnarray}
We see that, if $\xi \ne 0$, then $\lambda =i{\sqrt{k_2}} \xi$ is a pure imaginary eigenvalue of $A$.
If $\xi = 0$, then $\lambda =i{\sqrt{2k_1}}$ is a pure imaginary eigenvalue of $A$. Consequently (see  \cite{tesc}), the solution 
\eqref{e14} of \eqref{g6} doesn't converge to zero when times $t$ goes to infinity.
\end{proof} 

\subsection{Case 2: $(\tau_1 ,\tau_2 ,\tau_3)=(0,1,0)$}

We present, first, our exponential stability result for \eqref{g6}. 
\vskip0,1truecm
\begin{lemma}\label{lemma7432}
The result of Lemma \ref{lemma432} is satisfied in case \eqref{e102} with $(\tau_1 ,\tau_2 ,\tau_3)=(0,1,0)$.
\end{lemma}
\vskip0,1truecm
\begin{proof}
We addapt the arguments used in Subsection 3.2. We define $F_0$-$F_4$ and $F$ as in Subsection 3.2, and we get \eqref{F4g26+}, where
$F_6$ is defined by \eqref{7F6N1} with $k_4\xi\eta$ replaced by $ik_4 q$. Let $\lambda_6 >0$ and
\begin{eqnarray*}
{\tilde F} (\xi,\,t) &=& \xi^2 {F} (\xi,\,t)+ \lambda_6 \xi^2 Re\left(i\xi\widehat{\eta}\, \overline{\widehat{q}} \right)+\frac{1}{k_4}\left(\frac{k_1}{\gamma}(I_1 -I_2)+(k_1 \lambda_5 -\gamma\lambda_4 )\xi^2 \right)\xi^2  
Re\left(\widehat{v}\, \overline{\widehat{q}} \right) \\
& & -\frac{1}{k_4}\left(I_3 +\frac{k_3}{\gamma}I_1\right) \xi^2 
Re\left(i\xi\widehat{\phi}\, \overline{\widehat{q}} \right)+\frac{1}{k_4}(k_2\lambda_5 -\gamma\lambda_1)\xi^2  
Re\left(i\xi\widehat{z}\, \overline{\widehat{q}} \right).
\end{eqnarray*}
Multiplying \eqref{7equ1}, \eqref{7equ3}, \eqref{7equ4}, \eqref{7equ6} and \eqref{F4g26+} by 
\begin{equation*}
\lambda_6\xi^2, \quad\frac{1}{k_4}\left(\frac{k_1}{\gamma}(I_1 -I_2)+k_1 \lambda_5 -\gamma \right)\xi^2, \quad
\frac{-1}{k_4}\left(I_3 +\frac{k_3}{\gamma}I_1\right)\xi^2, \quad \frac{1}{k_4}(\gamma\lambda_1 -k_2\lambda_5 )\xi^2\quad\hbox{and}\quad \xi^2,
\end{equation*} 
respectively, adding the obtained equations and applying Young's inequality for the terms depending on $\widehat{q}$, we find, for any $\varepsilon_0 >0$,
\begin{equation} \label{7F43}
\frac{d}{dt}{\tilde F}(\xi,\,t) \leq C_{\varepsilon_0 ,\lambda_1 ,\cdots,\lambda_6} {\tilde f}(\xi) \vert\widehat{q}\vert^{2} 
- \left(\gamma\lambda_5 -\lambda_1 - 1-\varepsilon_0\right)\xi^4\vert\widehat{y}\vert^{2}
- \left(k_4\lambda_6 -\gamma\lambda_5 -\varepsilon_0\right)\xi^4\vert\widehat{\eta}\vert^{2}
\end{equation} 
\begin{equation*}
-\xi^{4} \,\left((k_1 -k_1 \lambda_2 -k_1 \lambda_{4} -\varepsilon_0 )\vert\widehat{v}\vert^{2} +(k_2 \lambda_1 -\varepsilon_0 )\vert\widehat{z}\vert^{2} +(\lambda_4 -\lambda_{3} -\varepsilon_0)\vert\widehat{\theta}\vert^{2} +(\lambda_2 -\varepsilon_0)\vert\widehat{u}\vert^{2} +(k_3 \lambda_3 -\varepsilon_0)\vert\widehat{\phi}\vert^{2}\right) ,
\end{equation*} 
where ${\tilde f}$ is defined in \eqref{4g270}. We choose $0<\lambda_{1}$, $0<\lambda_{2}<1$, $0<\lambda_{3}<\lambda_{4} <1 -\lambda_2$,  $\lambda_{5} >\frac{1}{\gamma}(\lambda_{1} +1)$, $\lambda_{6} >\frac{\gamma}{k_4}\lambda_{5}$ and 
\begin{equation*}
0<\varepsilon_{0} <\min\left\{k_1 -k_1 \lambda_2 -k_1 \lambda_4 , k_2 \lambda_1,\lambda_4 -\lambda_3 ,
\lambda_2 ,k_3\lambda_3, \gamma\lambda_5 -\lambda_1 -1, k_4 \lambda_{6} -\gamma \lambda_{5}\right\}. 
\end{equation*}
Thus, using the definition of $\widehat{E}$, \eqref{7F43} implies that, for some positive constant $c_1$,
\begin{equation}\label{74g27+} 
\frac{d}{dt}{\tilde F}(\xi,\,t) \leq -c_1 \xi^{4}\widehat{E} (\xi,\,t)+ C{\tilde f}(\xi)\vert\widehat{q}\vert^{2} .
\end{equation}
Therefore, we introduce the functional
\begin{equation}\label{74g27}
{L}(\xi,\,t) = \lambda\,\widehat{E}(\xi,\,t) +\frac{1}{{\tilde f}(\xi) }{\tilde F}(\xi,\,t)
\end{equation}
and we deduce that, using \eqref{Ep22} and \eqref{74g27+} 
\begin{eqnarray}\label{74g30}
\frac{d}{dt}{L}(\xi,\,t) \leq -c_1 f(\xi)\widehat{E}(\xi,\,t)-\left(k_5\,\lambda - C\right)\vert\widehat{q}\vert^{2} , 
\end{eqnarray}
where $f$ is defined in \eqref{422}. The proof can be ended as for Lemma \ref{lemma32}.
\end{proof}
\vskip0,1truecm
\begin{theorem}\label{theorem7441}
The result of Theorem \ref{theorem441} is staisfied in case \eqref{s22} with $(\tau_1 ,\tau_2 ,\tau_3 )=(0,1,0)$. 
\end{theorem}
\vskip0,1truecm
\begin{proof} The proof is identical to the one of Theorem \ref{theorem441}. 
\end{proof}

\subsection{Case 3: $(\tau_1 ,\tau_2 ,\tau_3)=(0,0,1)$}

In this case, we prove the same stability results for \eqref{g6} and \eqref{e11} that given in Subsection 4.2, and moreover, the proofs are very similar. 
\vskip0,1truecm
\begin{lemma}\label{lemma7532}
The result of Lemma \ref{lemma432} holds true also in case \eqref{e102} with $(\tau_1 ,\tau_2 ,\tau_3)=(0,0,1)$. 
\end{lemma}
\vskip0,1truecm
\begin{proof}
We define $F_0$-$F_4$ and $F$ as in Subsection 3.3 and we obtain \eqref{FF4g26+}, where
$F_6$ is defined in \eqref{7F6N} with $k_4\xi\eta$ replaced by $ik_4 q$. Let $\lambda_6 >0$ and
\begin{eqnarray*}
{\tilde F} (\xi,\,t) &=& \xi^2 {F} (\xi,\,t)+ \lambda_6 \xi^2 Re\left(i\xi\widehat{\eta}\, \overline{\widehat{q}} \right)+\frac{1}{k_4}\left(\frac{1}{\gamma}(k_1 I_3 -I_1)+(k_1 \lambda_5 -\gamma\lambda_4)\xi^2 \right)\xi^2  
Re\left(\widehat{v}\, \overline{\widehat{q}} \right) \\
&& +\frac{1}{k_4}\left(\gamma\lambda_3 -k_3\lambda_5\right) \xi^2 
Re\left(i\xi\widehat{\phi}\, \overline{\widehat{q}} \right)-\frac{1}{k_4}\left(I_2 +\frac{k_2}{\gamma}I_3\right)\xi^2  
Re\left(i\xi\widehat{z}\, \overline{\widehat{q}} \right).
\end{eqnarray*}
Multiplying \eqref{7equ1}, \eqref{7equ3}, \eqref{7equ4}, \eqref{7equ6} and \eqref{FF4g26+} by 
\begin{equation*}
\lambda_6\xi^2, \quad\frac{1}{k_4}\left(\frac{k_1 }{\gamma}(I_3 -I_1)-(k_1 \lambda_5 +\gamma\lambda_4)\xi^2 \right)\xi^2, 
\quad\frac{1}{k_4}\left(\gamma\lambda_3 -k_3\lambda_5\right)\xi^2, \quad -\frac{1}{k_4}\left(I_2 +\frac{k_2}{\gamma}I_3\right)\xi^2\quad\hbox{and}\quad\xi^2,
\end{equation*}
respectively, adding the obtained equations and applying Young's inequality for the terms depending on $\widehat{q}$, we find, for any $\varepsilon_0 >0$,
\begin{equation} \label{7FF43}
\frac{d}{dt}{\tilde F}(\xi,\,t) \leq C_{\varepsilon_0 ,\lambda_1 ,\cdots,\lambda_6} {\tilde f}(\xi) \vert\widehat{q}\vert^{2} -\xi^4 \left(\gamma\lambda_5 -\lambda_3 - \lambda_4 -\varepsilon_0\right)\vert\widehat{\theta}\vert^{2} - \xi^4 \left(k_4\lambda_6 -\gamma\lambda_5 -\varepsilon_0\right)\vert\widehat{\eta}\vert^{2}
\end{equation} 
\begin{equation*}
-\xi^{4} \,\left((k_1 \lambda_4 -k_1\lambda_{2} -k_1 -\varepsilon_0 )\vert\widehat{v}\vert^{2}+(k_2 \lambda_1 -\varepsilon_0 ) \vert\widehat{z}\vert^{2} +(1-\lambda_1 -\varepsilon_0 )\vert\widehat{y}\vert^{2} + (\lambda_2 -\varepsilon_0 ) \vert\widehat{u}\vert^{2} +(k_3 \lambda_3 -\varepsilon_0 )\vert\widehat{\phi}\vert^{2}\right) ,
\end{equation*}
${\tilde f}$ is defined in \eqref{4g270}.
We choose $0<\lambda_{3}$, $0<\lambda_{1}<1$, $\lambda_{4} >1$, $0<\lambda_{2}<\lambda_{4} -1$, 
$\lambda_{5} >\frac{1}{\gamma} (\lambda_3 +\lambda_4 )$, $\lambda_{6} >\frac{\gamma}{k_4} \lambda_5$ and 
\begin{equation*}
0<\varepsilon_{0} <\min\left\{k_3\lambda_3 , \lambda_2,
1-\lambda_1 ,k_2\lambda_1 , k_1\lambda_4 -k_1\lambda_2 -k_1, \gamma \lambda_5 -\lambda_3 -\lambda_4, k_4\lambda_{6} -\gamma\lambda_5\right\}. 
\end{equation*}
Then, using the definition of $\widehat{E}$, \eqref{7FF43} imlies \eqref{74g27+}, and then \eqref{74g30} holds true. 
Consequentely, the proof can be ended as for Lemma \ref{lemma7432}.
\end{proof}
\vskip0,1truecm
\begin{theorem}\label{theorem7541}
The stability result given in Theorem \ref{theorem7441} is satisfied when $(\tau_1 ,\tau_2 ,\tau_3)=(0,0,1)$.
\end{theorem}
\vskip0,1truecm
\begin{proof} 
The proof is identical to the one of Theorem \ref{theorem7441}.
\end{proof}

\section{Comments and issues}

1. The optimality of the obtained decay rates on $\Vert\partial_{x}^{j}U\Vert_{L^{2}(\mathbb{R})}$ is an interesting open question.
This question will be the focus of our attention in a future work.
\vskip0,1truecm
2. When $(\tau_1 , \tau_2 ,\tau_3)\in \{(0,1,0), (0,0,1)\}$ and $k_1 =k_2 =k_3$, the function $f$ tends to $1$ when $\xi$ goes to infinity, which avoid the regularity loss property; that is, \eqref{estim2} with $j=\ell=0$ gives the stability of \eqref{s21} and \eqref{s22} with a decay rate of $\Vert U\Vert_{L^{2}(\mathbb{R})}$ depending only on 
$\Vert U_{0} \Vert_{L^{1}(\mathbb{R})}$ and $\Vert U_{0} \Vert_{L^{2}(\mathbb{R})}$.  
However, in the other cases, $f$ tends to $0$ when $\xi$ goes to infinity, this means that the dissipation is very weak in the high frequency region, which imposes the regularity loss property in the estimates because \eqref{estim1} and \eqref{estim3} with $j=\ell=0$ imply only the boundedness of $\Vert U\Vert_{L^{2}(\mathbb{R})}$.
\vskip0,1truecm
3. The estimate \eqref{estim2} leads to a faster speed of convergence to zero of $\Vert\partial_{x}^{j}U\Vert_{L^{2}(\mathbb{R})}$ than the one guareented by \eqref{estim1} and \eqref{estim3}. This can be explained by the fact that the Cattaneo law generates a dissipation stronger than the one generated by the Fourier law. On the other hand, for both laws with 
$(\tau_1 , \tau_2 ,\tau_3)\in \{(0,1,0), (0,0,1)\}$, the situation is more favorable when $k_1 =k_2 =k_3$ than in the opposite case.
\vskip0,1truecm
4. From the mathematical point of view, one can take $\gamma\in\mathbb{R}^*$ in case \eqref{s21}, and $\gamma,\,k_4\in\mathbb{R}^*$ in case \eqref{s22} (instead of $\gamma,\,k_4 >0$). The unique needed modifications of proofs when $\gamma,\,k_4 <0$ are multiplying \eqref{equ1}, \eqref{equ4}, \eqref{equ7} and \eqref{7equ1} by $-1$, and using the obtained identities instead of \eqref{equ1}, \eqref{equ4}, \eqref{equ7} and \eqref{7equ1}.   
\vskip0,1truecm
5. The coupling terms 
\begin{equation}\label{coupling1}
\tau_j \gamma\eta_x\quad\hbox{and}\quad \gamma(\tau_1 \varphi_{xt} +\tau_2 \psi_{xt} +\tau_3 w_{xt})
\end{equation} 
in \eqref{s21} and \eqref{s22} are of order one with respect to $x$. Mathematicaly, these coupling terms can be replaced by (order zero with respect to $x$) 
\begin{equation}\label{coupling0}
\tau_j \gamma\eta\quad\hbox{and}\quad -\gamma(\tau_1 \varphi_{t} +\tau_2 \psi_{t} +\tau_3 w_{t}),
\end{equation}
respectively, with $\gamma\in\mathbb{R}^*$. In this case, the terms 
$i\tau_j \gamma\xi\widehat{\eta}$ and $i\gamma\xi(\tau_1 \widehat{u} +\tau_2 \widehat{y} +\tau_3 \widehat{\theta})$ in \eqref{e101} and \eqref{e102} are replaced by $\tau_j \gamma\widehat{\eta}$ and $-\gamma(\tau_1 \widehat{u} +\tau_2 \widehat{y} +\tau_3 \widehat{\theta})$, respectively. On the other hand, \eqref{11} holds true with   
\begin{equation*}
f(\xi) = \frac{\xi^{8}}{1 + \xi^{2} +\xi^{4} +\xi^{6} +\xi^{8} +\xi^{10}}
\end{equation*} 
instead of \eqref{22}, and 
\begin{align*}
f(\xi) = \begin{cases}
\frac{\xi^{6}}{1 +\xi^{2} +\xi^{4} +\xi^{6}}\quad &\hbox{if}\,\, k_1 =k_2 =k_3 , \\
\frac{\xi^{6}}{1 +\xi^{2} +\xi^4 +\xi^6 +\xi^{8} +\xi^{10}} \quad &\hbox{if not}
\end{cases}
\end{align*}
instead of \eqref{422}, and so we get the stability estimates 
\begin{equation*}
\Vert\partial_{x}^{j}U\Vert_{L^{2}(\mathbb{R})} \leq c_0 \,(1 + t)^{-1/16 - j/8}\,\Vert U_{0} \Vert_{L^{1}(\mathbb{R})} + c_0 \,(1 + t)^{-\ell/2}\,\Vert\partial_{x}^{j+\ell}U_{0} \Vert_{L^{2}(\mathbb{R})},\quad\forall t\in \mathbb{R}_+ 
\end{equation*}  
instead of \eqref{33},
\begin{equation*}
\Vert\partial_{x}^{j}U\Vert_{L^{2}(\mathbb{R})} \leq c_0 \,(1 + t)^{-1/12 - j/6}\,\Vert U_{0} \Vert_{L^{1}(\mathbb{R})} + c_0 e^{-{\tilde c}_0 t}\,\Vert\partial_{x}^{j+\ell}U_{0} \Vert_{L^{2}(\mathbb{R})} \quad\hbox{if}\,\,k_1 =k_2 =k_3  
\end{equation*} 
instead of \eqref{33case20}, and 
\begin{equation*}
\Vert\partial_{x}^{j}U\Vert_{L^{2}(\mathbb{R})} \leq c_0 \,(1 + t)^{-1/12 - j/6}\,\Vert U_{0} \Vert_{L^{1}(\mathbb{R})} + c_0 \,(1 + t)^{-\ell/4}\,\Vert\partial_{x}^{j+\ell}U_{0} \Vert_{L^{2}(\mathbb{R})} \quad\hbox{if not} 
\end{equation*}    
instead of \eqref{33case2}. These stability estimates show that the decay rates in case \eqref{coupling0} is smaller than the ones obtained in case \eqref{coupling1}. Moreover, the non stability result when $k_2 =k_3$ and $(\tau_1 ,\tau_2 ,\tau_3 )=(1,0,0)$ is still valid using the same arguments of proof, since we get \eqref{detlambdaIA} and \eqref{detlambdaIA2} with $\gamma^2$ instead of $\gamma^2 \xi^2$.


\begin{thebibliography}{1}
\bibitem{agg} M. S. Alves, P. Gamboa, G. C. Gorain, A. Rambaud and O. Vera, Asymptotic behavior of a flexible structure with Cattaneo type of thermal effect, Indagationes Mathematicae, 27 (2016), 821-834.

\bibitem{beim} C. F. Beards and I. M. A. Imam, The damping of plate vibration by interfacial slip between layers, Int. J. Mach. Tool. Des. Res., 18 (1978), 131-137. 

\bibitem{clx} X. G. Cao, D. Y. Liu and G. Q. Xu, Easy test for stability of laminated beams with structural
damping and boundary feedback controls, J. Dynamical Control Syst., 13 (2007), 313-336.

\bibitem{cdflr} M. M. Cavalcanti, V. N. Domingos Cavalcanti, F. A. Falcao Nascimento, I. Lasiecka and J. H. Rodrigues, Uniform decay rates for the energy of Timoshenko system with the arbitrary speeds of propagation and localized nonlinear damping, Z. Angew. Math. Phys., 65 (2014), 1189-1206.

\bibitem{s6} L. Djouamai and B. Said-Houari, A new stability number of the Bresse-Cattaneo system, 
Math. Meth. Appl. Sci., 41 (2018), 2827-2847.

\bibitem{4} L. H. Fatori, R. N. Monteiro and H. D. Fern\a'{a}ndez Sare, The Timoshenko system with history and Cattaneo
law, Applied Mathematics and Computation, 228 (2014), 128-140.

\bibitem{s1} T. E. Ghoul, M. Khenissi and B. Said-Houari, On the stability of the Bresse system with frictional damping, J. Math. Anal. Appl., 455 (2017), 1870-1898.

\bibitem{gues2} A. Guesmia, Asymptotic stability of Bresse system with one infinite memory in the longitudinal displacements, Medi. J. Math., 14 (2017), 19 pages.

\bibitem{gues4} A. Guesmia, Non-exponential and polynomial stability results of a Bresse system with one infinite memory in the vertical displacement, Nonauton. Dyn. Syst., 4 (2017), 78-97.

\bibitem{gues8} A. Guesmia, Well-posedness and stability results for laminated Timoshenko beams with interfacial slip and infinite memory, IMA J. Math. Cont. Info., 37 (2020), 300-350.

\bibitem{gms} A. Guesmia, S. Messaoudi and A. Soufyane, On the stabilization for a linear Timoshenko system with infinite history and applications to the coupled Timoshenko-heat systems, Elec. J. Diff. Equa., 2012 (2012), 1-45.

\bibitem{hans} S. W. Hansen, In control and estimation of distributed parameter systems: Non-linear
phenomena, International Series of Numerical Analysis, 118 (1994), 143-170.

\bibitem{hasp} S. W. Hansen and R. Spies, Structural damping in a laminated beams due to interfacial slip,
J. Sound Vibration, 204 (1997), 183-202.

\bibitem{8} K. Ide, K. Haramoto and S. Kawashima, Decay property of regularity-loss type for dissipative Timoshenko
system, Math. Mod. Meth. Appl. Sci., 18 (2008), 647-667.

\bibitem{s4} M. Khader and B. Said-Houari, Decay rate of solutions to Timoshenko system with past history in unbounded domains, 
Appl. Math. Optim., 75 (2017), 403-428.

\bibitem{s5} M. Khader and B. Said-Houari, Optimal decay rate of solutions to Timoshenko system with past history in unbounded domains, Z. Anal. Anwend, 37 (2018), 435-459.

\bibitem{lizh} W. Liu and W. Zhao, Exponential and polynomial decay for a laminated beam with Fourier's type heat conduction, Preprints 2017, 2017020058, doi: 10.20944/preprints201702.0058.v1.

\bibitem{lota1} A. Lo and N-E Tatar, Stabilization of laminated beams with interfacial slip, Elec. J. Diff. Equa., 2015 (2015), 1-14.

\bibitem{lota2} A. Lo and N. E. Tatar, Uniform stability of a laminated beam with structural memory, Qual. Theory Dyn.
Syst., 15 (2016), 517-540.

\bibitem{lota3} A. Lo and N. E. Tatar, Exponential stabilization of a structure with interfacial slip, Discrete Contin. Dyn.
Syst., 36 (2016), 6285-6306.

\bibitem{rapo} C. A. Raposo, Exponential stability for a structure with interfacial slip and
frictional damping, Appl. Math. Lett., 53 (2016), 85-91.

\bibitem{rvma} C. A. Raposo, O. V. Villagr\'an, J. E. Mu\~noz Rivera and M. S. Alves, Hybrid laminated Timoshenko beam, J. Math. Phys., 58 (2017), 11 pages. 

\bibitem{s3} B. Said-Houari and R. Racke, Decay rates and global existence for semilinear dissipative Timoshenko systems, 
Quart. Appl. Math., 71 (2013), 229-266.

\bibitem{6} B. Said-Houari and R. Rahali, Asymptotic behavior of the Cauchy problem of the Timoshenko system
in thermoelsaticity of type III, Evolution Equations and Control Theory, 2(2013), 423-440.

\bibitem{s2} B. Said-Houari and A. Soufyane, The effect of frictional damping terms on the decay rate of the Bresse system, Evol. Equa. Cont. Theory, 3 (2014), 713-738.

\bibitem{10} M. L. Santos, D. S. Almeida and J. E. Mu\~noz Rivera, The stability number of the Timoshenko system
with second sound, J. Diff. Equa., 253 (2012), 2715-2733.

\bibitem{7} A. Soufyane and B. Said-Houari, The effect of the wave speeds and the frictional damping terms on the
decay rate of the Bresse system, Evolution Equations and Control Theory, 3 (2014), 713-738.

\bibitem{tesc} G. Teschl, Ordinary differential equations and dynamical systems, American Mathematical
Soc., 140 (2012), ISBN 978-0-8218-8328-0.

\bibitem{wxy}  J. M. Wang, G. Q. Xu and S. P. Yung, Exponential stabilization of laminated beams with structural damping and boundary feedback controls, SIAM J. Control Optim., 44 (2005), 1575-1597.
\end{thebibliography}
\end{document}